\documentclass[largeformat]{interact}

\usepackage{enumitem}
\usepackage{amsmath}
\usepackage{mathtools}
\newcommand{\ie}{i.\,e.~}
\newcommand{\eg}{e.\,g.~}
\newcommand{\Span}[1]{\mathrm{span}\{#1\}}
\newcommand{\Rank}[1]{\mathrm{rank}(#1)}
\newcommand{\Dim}[1]{\mathrm{dim}(#1)}
\newcommand{\C}[1]{\mathcal{C}(#1)}
\newcommand{\D}{\mathrm{d}}
\newcommand{\Lie}{\mathrm{L}}
\newcommand{\Ad}{\mathrm{ad}}

\newcommand{\Mod}{\mathrm{~mod~}}

\newcounter{fact}[subsection]
\newenvironment{fact}[1][]{\refstepcounter{fact}\par\medskip
	\noindent\textit{Fact~\thefact. #1} \itshape}{\medskip\newline}
\newenvironment{factProof}[1][]{\par\medskip
	\noindent\textit{Proof of Fact #1.} \rmfamily}{\medskip\newline}

\newcounter{prop}[subsection]
\newenvironment{prop}[1][]{\refstepcounter{prop}\par\medskip
	\noindent\textit{Proposition~\theprop. #1} \itshape}{\medskip\newline}
\newenvironment{propProof}[1][]{\par\medskip
	\noindent\textit{Proof of Proposition #1.} \rmfamily}{\medskip\newline}

\usepackage{epstopdf}
\usepackage[caption=false]{subfig}

\usepackage[natbibapa,nodoi]{apacite}
\theoremstyle{plain}
\newtheorem{theorem}{Theorem}[section]
\newtheorem{lemma}[theorem]{Lemma}

\theoremstyle{definition}

\theoremstyle{remark}
\newtheorem{remark}{Remark}

\begin{document}

\articletype{ARTICLE TEMPLATE}

\title{A Structurally Flat Triangular Form Based on the Extended Chained Form}

\author{
	\name{Conrad Gst\"ottner\textsuperscript{a}\thanks{CONTACT Conrad Gst\"ottner. Email: conrad.gstoettner@jku.at\newline The first author and the second author have been supported by the Austrian Science Fund (FWF) under grant number P 32151 and P 29964.}, Bernd Kolar\textsuperscript{a} and Markus Sch\"oberl\textsuperscript{a}}
	\affil{\textsuperscript{a}Institute of Automatic Control and Control Systems Technology, Johannes Kepler University, Linz, Austria;
	}
}

\maketitle

\begin{abstract}
	In this paper, we present a structurally flat triangular form which is based on the extended chained form. We provide a complete geometric characterization of the proposed triangular form in terms of necessary and sufficient conditions for an affine input system with two inputs to be static feedback equivalent to this triangular form. This yields a sufficient condition for an affine input system to be flat.
\end{abstract}


	\section{Introduction}
		The concept of flatness was introduced in control theory by Fliess, L\'evine, Martin and Rouchon, see \eg \cite{FliessLevineMartinRouchon:1992,FliessLevineMartinRouchon:1995}, and has attracted a lot of interest in the control systems theory community. The flatness property allows an elegant systematic solution of feed-forward and feedback problems, see \eg \cite{FliessLevineMartinRouchon:1995}. Roughly speaking, a nonlinear control system
		\begin{align}\label{eq:sys}
			\begin{aligned}
				\dot{x}&=f(x,u)
			\end{aligned}
		\end{align}
		with $\Dim{x}=n$ states and $\Dim{u}=m$ inputs is flat, if there exist $m$ differentially independent functions $y^j=\varphi^j(x,u,u_1,\ldots,u_q)$, $u_k$ denoting the $k$-th time derivative of $u$, such that $x$ and $u$ can be parameterized by $y$ and its time derivatives. Up to now, there do not exist verifiable necessary and sufficient conditions for testing a system of the form \eqref{eq:sys} for flatness, only for certain subclasses of systems, the flatness problem has been solved. Recent research in the field of flatness can be found in \eg \cite{SchoberlRiegerSchlacher:2010}, \cite{SchlacherSchoberl:2013}, \cite{LiXuSuChu:2013}, \cite{SchoberlSchlacher:2014}, \cite{KolarSchoberlSchlacher:2015}, \cite{NicolauRespondek:2017}.\\\\
		Structurally flat triangular forms are of special interest in the problem of deriving flat outputs for nonlinear control systems. In \cite{BououdenBoutatZhengBarbotKratz:2011}, a structurally flat triangular form for a class of $0$-flat systems is proposed and geometric necessary and sufficient conditions for the existence of a transformation of a nonlinear control system into this triangular form are provided. A structurally flat implicit triangular form for $1$-flat systems, together with a constructive scheme for transforming a system into the proposed triangular form, can be found in \cite{SchoberlSchlacher:2014}. A complete solution of the flatness problem of systems that become static feedback linearizable after a one-fold prolongation of a suitably chosen control is presented in \cite{NicolauRespondek:2017}. Normal forms for this class of systems can be found in \cite{NicolauRespondek:2019}. Another class of systems for which the flatness problem has been solved are two-input driftless systems, see \cite{MartinRouchon:1994}. Flat two-input driftless systems are static feedback equivalent to a structurally flat triangular form, referred to as chained form. In \cite{LiXuSuChu:2013} an extension of the chained form for systems with drift, the so called extended chained form, is considered. Geometric necessary and sufficient conditions for a two-input affine input system (AI-system) to be static feedback equivalent to this extended chained form can be found in \cite{SilveiraPereiraRouchon:2015}. Conditions for the case with $m\geq 2$ inputs are provided in \cite{Nicolau:2014}, \cite{NicolauLiRespondek:2014} and \cite{LiNicolauRespondek:2016}.\\\\	
		In \cite{GstottnerKolarSchoberl:2020}, a triangular form which generalizes the extended chained form is considered and necessary and sufficient conditions for an two-input AI-system to be static feedback equivalent to this triangular form are provided. The proposed triangular form generalizes the extended chained form by augmenting it with two subsystems in Brunovsky normal form. To be precise, two equally lengthened integrator chains are attached to the inputs of a subsystem in extended chained form and furthermore, the top variables (flat outputs) of this subsystem in extended chained form act as inputs for two arbitrary lengthened integrator chains. The structurally flat triangular form obtained this way contains the (extended) chained form as a special case. In this contribution, we further develop the ideas presented in \cite{GstottnerKolarSchoberl:2020}. We again augment the extended chained form with integrator chains, but here, the integrator chains attached to the inputs of the subsystem in (extended) chained form differ in length by one integrator. As a consequence, the (extended) chained form is not contained as a special case. It turns out that a broad variety of practical and academic examples is static feedback equivalent to this triangular form. Among others, \eg the planar VTOL aircraft, also considered in \eg \cite{FliessLevineMartinRouchon:1999} and \cite{SchoberlRiegerSchlacher:2010} and the model of a gantry crane, considered in \eg \cite{FliessLevineMartinRouchon:1995}. These systems cannot be handled with the triangular form presented in \cite{GstottnerKolarSchoberl:2020}. We again provide necessary and sufficient conditions for an AI-system to be static feedback equivalent to this triangular form. This again provides a sufficient condition for an AI-system to be flat. In contrast to \cite{GstottnerKolarSchoberl:2020}, where proofs are only sketched, detailed proofs are provided in this contribution.
	\section{Notation}\label{se:not}
		Let $\mathcal{X}$ be an $n$-dimensional smooth manifold, equipped with local coordinates $x^i$, $i=1,\ldots,n$. Its tangent bundle and cotangent bundle are denoted by $(\mathcal{T}(\mathcal{X}),\tau_\mathcal{X},\mathcal{X})$ and $(\mathcal{T}^\ast(\mathcal{X}),\tau^\ast_\mathcal{X},\mathcal{X})$. For these bundles we have the induced local coordinates $(x^i,\dot{x}^i)$ and $(x^i,\dot{x}_i)$ with respect to the bases $\{\partial_{x^i}\}$ and $\{\D x^i\}$, respectively. Throughout, the Einstein summation convention is used. The exterior derivative of a $p$-form $\omega$ is denoted by $\D\omega$. By $\Lie_v^k\varphi$ we denote the $k$-fold Lie derivative of a function $\varphi$ along a vector field $v$. Let $v$ and $w$ be two vector fields. Their Lie bracket is denoted by $[v,w]$, for the repeated application of the Lie bracket, we use the common notation $\Ad_v^kw=[v,\Ad_v^{k-1}w]$, $k\geq 1$ and $\Ad_v^0w=w$. Let furthermore $D_1$ and $D_2$ be two distributions. By $[v,D_1]$ we denote the distribution spanned by the Lie bracket of $v$ with all basis vector fields of $D_1$, and by $[D_1,D_2]$ the distribution spanned by the Lie brackets of all possible pairs of basis vector fields of $D_1$ and $D_2$. The $i$-th derived flag of a distribution $D$ is denoted by $D^{(i)}$ and defined by $D^{(0)}=D$ and $D^{(i+1)}=D^{(i)}+[D^{(i)},D^{(i)}]$ for $i\geq 0$. The $i$-th Lie flag of a distribution $D$ is denoted by $D_{(i)}$ and defined by $D_{(0)}=D$ and $D_{(i+1)}=D_{(i)}+[D,D_{(i)}]$ for $i\geq 0$. The involutive closure of $D$ is denoted by $\overline{D}$, it is the smallest involutive distribution which contains $D$. It can be determined via the derived flag. We denote the Cauchy characteristic distribution of $D$ by $\C{D}$. It is spanned by all vector fields $c$ which belong to $D$ and satisfy $[c,D]\subset D$. Cauchy characteristic distributions are always involutive. They allow us to find a basis for a distribution which is independent of certain coordinates. Since $\C{D}$ is involutive, it can be straightened out such that $\C{D}=\Span{\partial_{x^1},\ldots,\partial_{x^{n_c}}}$, with $n_c=\Dim{\C{D}}$. From $[\C{D},D]\subset D$, it follows that in these coordinates, a basis for $D$ which does not depend on the coordinates $(x^1,\ldots,x^{n_c})$ can be constructed. Consider an AI-system with $m$-inputs
		\begin{align}\label{eq:aisystemIntro}
			\begin{aligned}
				\dot{x}&=a(x)+b_j(x)u^j\,,&j&=1,\ldots,m\,.
			\end{aligned}
		\end{align}
		Geometrically, such a system is represented by the drift vector field $a=a^i(x)\partial_{x^i}$ and the input vector fields $b_j=b_j^i(x)\partial_{x^i}$, $j=1,\ldots,m$, $i=1,\ldots,n$ on the state manifold $\mathcal{X}$. Throughout, we assume that all vector fields and functions we deal with are smooth. We call two AI-systems static feedback equivalent, if they are equivalent via a diffeomorphism $\tilde{x}=\Phi(x)$ on the state space and an invertible feedback transformation $\tilde{u}^j=g^j(x)+m^j_k(x)u^k$. The equivalent system reads
		\begin{align*}
			\begin{aligned}
				\dot{\tilde{x}}&=\tilde{a}(\tilde{x})+\tilde{b}_j(\tilde{x})\tilde{u}^j\,,&j&=1,\ldots,m\,.
			\end{aligned}
		\end{align*}
	\section{Known Results}\label{se:knownResults}
		In this section, we summarize some known results from the literature which are of particular importance for characterizing our triangular form. Throughout, we assume all distributions to have locally constant dimension, we consider generic points only. In particular, we call a system static feedback equivalent to a certain normal form, even though the transformation into this form may exhibit singularities. Consider again an $m$-input AI system \eqref{eq:aisystemIntro}. Such a system is called static feedback linearizable, if it is static feedback equivalent to a linear controllable system, in particular to the Brunovsky normal form. The static feedback linearization problem has been solved in \cite{JakubczykRespondek:1980} and \cite{HuntSu:1981}. The geometric necessary and sufficient conditions read as follows. For \eqref{eq:aisystemIntro}, we define the distributions $D_{i+1}=D_i+[a,D_i]$, $i\geq 1$, where $D_1=\Span{b_1,\ldots,b_m}$.
		\begin{theorem}\label{thm:exactLin}
			The $m$-input AI-system \eqref{eq:aisystemIntro} is static feedback linearizable if and only if all the distributions $D_i$, $i\geq 1$ are involutive and $D_{n-1}=\mathcal{T}(\mathcal{X})$.
		\end{theorem}\noindent
		In \cite{MartinRouchon:1994} it is shown that a two-input driftless system of the form
		\begin{align}\label{eq:driftless}
			\begin{aligned}
				\dot{x}&=b_1(x)u^1+b_2(x)u^2
			\end{aligned}
		\end{align}
		is flat, if and only if it is static feedback equivalent to the structurally flat triangular form
		\begin{align}\label{eq:chainedform}
			\begin{aligned}
			\dot{x}^1=u^2,~\dot{x}^2=x^3u^2,~\cdots,~\dot{x}^{n-1}=x^nu^2,~\dot{x}^n=u^1\,,
			\end{aligned}
		\end{align}
		referred to as chained form. The input vector fields of a system in chained form read
		\begin{align}\label{eq:chaindVectorfields}
			\begin{aligned}
				b_1&=\partial_{x^n}\,,&b_2&=\partial_{x^1}+x^3\partial_{x^2}+\ldots+x^n\partial_{x^{n-1}}\,.
			\end{aligned}
		\end{align}
		The geometric necessary and sufficient conditions for a driftless system \eqref{eq:driftless} to be static feedback equivalent to the chained form \eqref{eq:chainedform} are summarized in the following theorem.
		\begin{theorem}\label{thm:chainedform}
			The driftless system \eqref{eq:driftless} is static feedback equivalent to the chained form \eqref{eq:chainedform} if and only if $D=\Span{b_1,b_2}$ satisfies $\Dim{D^{(i)}}=2+i$, $i=0,\ldots,n-2$.
		\end{theorem}\noindent
		In \cite{Murray:1994} it is shown that locally around a point of the state space at which the additional regularity condition $\Dim{D_{(i)}}=2+i$, $i=0,\ldots,n-2$ on the Lie flag of $D$ holds, the transformation into chained form does not exhibit singularities. A system in chained form is flat with the pair of top variables $(x^1,x^2)$ forming a possible flat output. For a comprehensive analysis of the flatness of systems static feedback equivalent to the chained form, a characterization of all their $x$-flat outputs and their singularities, we refer to \cite{LiRespondek:2012}. The structurally flat triangular form
		\begin{align}\label{eq:extendedChainedform}
			\begin{aligned}
				\dot{x}^1&=u^2\\
				\dot{x}^2&=x^3u^2+a^2(x^1,x^2,x^3)\\
				\dot{x}^3&=x^4u^2+a^3(x^1,x^2,x^3,x^4)\\
				&\vdotswithin{=}\\
				\dot{x}^{n-1}&=x^nu^2+a^{n-1}(x^1,\ldots,x^n)\\
				\dot{x}^n&=u^1\,,
			\end{aligned}
		\end{align}
		referred to as extended chained form, was first considered in \cite{LiXuSuChu:2013}. In \cite{Nicolau:2014} and \cite{SilveiraPereiraRouchon:2015} geometric necessary and sufficient conditions for an AI-system with two inputs to be static feedback equivalent to \eqref{eq:extendedChainedform} are provided. Those are summarized in the following theorem.
		\begin{theorem}\label{thm:extendedChainedform}
			An AI-system \eqref{eq:aisystemIntro} with two inputs ($m=2$) is static feedback equivalent to the extended chained form \eqref{eq:extendedChainedform} if and only if
			\begin{enumerate}[label=\arabic*)]
				\item $D=\Span{b_1,b_2}$ satisfies $\Dim{D^{(i)}}=\Dim{D_{(i)}}=2+i$, $i=0,\ldots,n-2$.
				\item The drift of the system meets the compatibility condition
				\begin{align}\label{eq:compatibility}
					\begin{aligned}
						[a,\mathcal{C}(D^{(i)})]&\subset D^{(i)}\,,&i&=1,\ldots,n-3\,.
					\end{aligned}
				\end{align}
			\end{enumerate}
		\end{theorem}\noindent
		These conditions can be interpreted as follows. The first condition assures that the drifless system obtained by setting $a(x)=0$ is static feedback equivalent to the chained form \eqref{eq:chainedform}, the condition $\Dim{D_{(i)}}=2+i$, $i=0,\ldots,n-2$ on the Lie flag of $D$ is again a regularity condition. The second condition assures that the drift is compatible with the chained form, \ie that in coordinates in which the input vector fields are in chained form, the drift takes the desired triangular structure. Thus, a flat output of the driftless system obtained by setting $a(x)=0$, is also a flat output of \eqref{eq:extendedChainedform}. For a comprehensive analysis of the flatness of systems static feedback equivalent to the extended chained form, a characterization of their flat outputs and their singularities, we refer to \cite{Nicolau:2014} or \cite{LiNicolauRespondek:2016}. 
		\begin{remark}\label{rem:linearization}
			Note that a system in chained form \eqref{eq:chainedform} becomes static feedback linearizable by $(n-2)$-fold prolonging the input $u^2$. The same holds for a system in extended chained form \eqref{eq:extendedChainedform}, see also \cite{Nicolau:2014}.
		\end{remark}
	\vspace{-1em}
	\section{A Structurally Flat Triangular Form Based on the Extended Chained Form}\label{se:tf}
		In the following, we consider the structurally flat triangular form
		\begin{align}\label{eq:triangularFormDifferentLengths}
			\begin{aligned}
				\dot{x}_1&=f_1(x_1,x_2^1,x_2^2)\\
				\dot{x}_2&=f_2(x_1,x_2,x_{3,1}^1,x_{3,2}^1)\\
				\dot{x}_3&=f_3(x_3,u^1,u^2)\,,
			\end{aligned}
		\end{align}
		with the $x_1$-subsystem being in Brunovsky normal form
		\begin{align}\label{eq:subsys1}
			f_1:\quad\begin{aligned}
				\dot{x}_{1,1}^1&=x_{1,1}^2&\dot{x}_{1,2}^1&=x_{1,2}^2\\
				\dot{x}_{1,1}^2&=x_{1,1}^3&\dot{x}_{1,2}^2&=x_{1,2}^3\\
				&\vdotswithin{=}&&\vdotswithin{=}\\
				\dot{x}_{1,1}^{n_{1,1}}&=x_2^1&\dot{x}_{1,2}^{n_{1,2}}&=x_2^2\,,
			\end{aligned}
		\end{align}
		the $x_2$-subsystem being essentially in extended chained form
		\begin{align}\label{eq:subsys2}
			f_2:\quad\begin{aligned}
				\dot{x}_2^1&=x_{3,2}^1\\
				\dot{x}_2^2&=x_2^3x_{3,2}^1+a_2^2(x_1,x_2^1,\ldots,x_2^3)\\
				&\vdotswithin{=}\\
				\dot{x}_2^{n_2-1}&=x_2^{n_2}x_{3,2}^1+a_2^{n_2-1}(x_1,x_2)\\
				\dot{x}_2^{n_2}&=x_{3,1}^1+g(x_1,x_2)x_{3,2}^1
			\end{aligned}
		\end{align}
		and the $x_3$-subsystem again being in Brunovsky normal form
		\begin{align}\label{eq:subsys3}
			f_3:\quad\begin{aligned}
				\dot{x}_{3,1}^1&=x_{3,1}^2&\dot{x}_{3,2}^1&=x_{3,2}^2\\
				\dot{x}_{3,1}^2&=x_{3,1}^3&\dot{x}_{3,2}^2&=x_{3,2}^3\\
				&\vdotswithin{=}&&\vdotswithin{=}\\
				\dot{x}_{3,1}^{n_3-1}&=x_{3,1}^{n_3}&\dot{x}_{3,2}^{n_3-1}&=u^2\\
				\dot{x}_{3,1}^{n_3}&=u^1\,.
			\end{aligned}
		\end{align}
		The triangular form \eqref{eq:triangularFormDifferentLengths} consists of three subsystems. The $x_1$-subsystem is in Brunovsky normal form, it consists of two integrator chains of arbitrary lengths $n_{1,1}\geq 0$ and $n_{1,2}\geq 0$. In total, it consists of $n_1=n_{1,1}+n_{1,2}\geq 0$ states. The $x_2$-subsystem is essentially in extended chained form (we assume $n_2\geq 3$) and the top variables $x_2^1$ and $x_2^2$ of this subsystem act as inputs for the $x_1$-subsystem. 
		\begin{remark}\label{rem:essentiallyExtendedChained}
			The $x_2$-subsystem differs from the extended chained form in two minor ways. Firstly, the functions $a_2^i$, $i=2,\ldots,n_2-1$, which represent the drift of the $x_2$-subsystem may also depend on the stats $x_1$. Secondly, in the last equation of the $x_2$-subsystem, besides $x_{3,1}^1$, there may also occur the term $g(x_1,x_2)x_{3,2}^1$. Nevertheless, the $x_2$-subsystem has analogous structural properties as a system in extended chained form.
		\end{remark}\noindent
		The $x_3$-subsystem is again in Brunovsky normal form, it consists of two integrator chains which differ in length by one integrator. The top variables $x_{3,1}^1$ and $x_{3,2}^1$ act as inputs for the $x_2$-subsystem\footnote{Note that the top variable $x_{3,1}^1$ of the longer integrator chain corresponds to the input $u^1$ in \eqref{eq:extendedChainedform}, \ie the input which only occurs in the very last equation of the extended chained form. This is crucial, if $x_{3,1}^1$ and $x_{3,2}^1$ would be swapped, the system would be static feedback equivalent to a system of the form \eqref{eq:triangularFormDifferentLengths} with equally lengthened integrator chains in the $x_3$-subsystem, \ie it would be static feedback equivalent to the triangular form proposed in \cite{GstottnerKolarSchoberl:2020}. In this case, the last equation of the $x_2$-subsystem would belong to the shorter integrator chain, which would compensate the length difference.\label{fn:wrongInput}} (we assume $n_3\geq 1$, for $n_3=1$, the $x_3$-subsystem only consists of a single integrator, namely $\dot{x}_{3,1}^1=u^1$, and in the $x_2$-subsystem, $x_{3,2}^1$ is replaced by $u^2$). In conclusion, the $x_3$-subsystem and the $x_2$-subsystem form an endogenous dynamic feedback for the $x_1$-subsystem. The $x_3$-subsystem in turn is an endogenous dynamic feedback for the $x_2$-subsystem. The total number of states of \eqref{eq:triangularFormDifferentLengths} is given by $n=n_1+n_2+2n_3-1$.
		\begin{remark}\label{rem:dimensions}
			In conclusion, the restriction on the dimensions of the subsystems in \eqref{eq:triangularFormDifferentLengths} are $n_2\geq 3$ and $n_3\geq 1$. However, it turns out that a system of the form \eqref{eq:triangularFormDifferentLengths} with $n_2=3$ and $n_1=0$ meets the conditions of Theorem \ref{thm:exactLin} and thus, it is static feedback linearizable. Therefore, $n_2=3$ only makes sense if $n_1\geq 1$.
		\end{remark}
		\begin{remark}
			A system of the form \eqref{eq:triangularFormDifferentLengths} becomes static feedback linearizable after an $(n_2-1)$-fold prolongation of $u^2$ (one prolongation accounts for the differing lengths of the integrator chains in the $x_3$-subsystem, the remaining $(n_2-2)$ prolongations correspond to those in Remark \ref{rem:linearization}). In particular, a system of the form \eqref{eq:triangularFormDifferentLengths} with $n_2=3$ becomes static feedback linearizable after a two-fold prolongation of $u^2$. A geometric characterization of systems that become static feedback linearizable after a two-fold prolongation of a
			suitably chosen control can be found in \cite{NicolauRespondek:2016-2}. However, due to Assumption 2 in \cite{NicolauRespondek:2016-2}, the geometric necessary and sufficient conditions for linearizability via a two-fold prolongation provided therein do not apply on a system of the form \eqref{eq:triangularFormDifferentLengths} if $n_{1,1}\leq 1$ and $n_{1,2}\leq 1$. The special case $n_2=3$ and $n_{1,1}\geq 2$ or $n_{1,2}\geq 2$ is indeed fully covered by \cite{NicolauRespondek:2016-2}. Our geometric characterization of \eqref{eq:triangularFormDifferentLengths} provided in the following section is not subject to any restrictions on $n_1$, $n_2$ and $n_3$ (except those in Remark \ref{rem:dimensions}).
		\end{remark}\noindent
		As a motivating example, consider the planar VTOL aircraft, also treated \eg in \cite{FliessLevineMartinRouchon:1999}, \cite{SchoberlRiegerSchlacher:2010} or \cite{SchoberlSchlacher:2011}, and given by
		\begin{align}\label{eq:vtol}
			\begin{aligned}
				\dot{x}&=v_x\\
				\dot{z}&=v_z\\
				\dot{\theta}&=\omega
			\end{aligned}&\qquad
			\begin{aligned}
				\dot{v}_x&=\epsilon\cos(\theta)u^2-\sin(\theta)u^1\\
				\dot{v}_z&=\cos(\theta)u^1+\epsilon\sin(\theta)u^2-1\\
				\dot{\omega}&=u^2\,.
			\end{aligned}
		\end{align}
		This system is not static feedback linearizable, but it is known to be flat. It is not static feedback equivalent to the triangular form proposed in \cite{GstottnerKolarSchoberl:2020}, but it is static feedback equivalent to the triangular form \eqref{eq:triangularFormDifferentLengths}. In Section \ref{se:ex}, we will systematically derive a state and input transformation, which brings \eqref{eq:vtol} into the form
		\begin{align}\label{eq:vtolTriangular}
			\begin{aligned}
				&f_1:\quad\begin{aligned}
					\dot{x}_{1,1}^1&=x_2^1\\
					\dot{x}_{1,2}^1&=x_2^2
				\end{aligned}&&\qquad f_2:\quad\begin{aligned}
					\dot{x}_2^1&=\tilde{u}^2\\
					\dot{x}_2^2&=x_2^3\tilde{u}^2+x_2^3\\
					\dot{x}_2^3&=x_{3,1}^1
				\end{aligned}&&\qquad f_3:\quad\begin{aligned}
					\dot{x}_{3,1}^1&=\tilde{u}^1\,,
				\end{aligned}
			\end{aligned}
		\end{align}
		which is of the form \eqref{eq:triangularFormDifferentLengths} with $n_{1,1}=n_{1,2}=1$, $n_2=3$ and $n_3=1$.
		\begin{remark}
			In \cite{NicolauRespondek:2020}, normal forms for systems that become static
			feedback linearizable after a two-fold prolongation of a suitably chosen control are presented. Therein, based on a given suitable flat output a representation of the VTOL analogous to \eqref{eq:vtolTriangular} is derived. Note however that the geometric necessary and sufficient conditions for linearizability via a two-fold prolongation provided in \cite{NicolauRespondek:2016-2} do not apply on the VTOL.
		\end{remark}\noindent
		The triangular form \eqref{eq:triangularFormDifferentLengths} is similar to the triangular form presented in \cite{GstottnerKolarSchoberl:2020}. The difference between the triangular form considered here and the triangular form in \cite{GstottnerKolarSchoberl:2020} is that in \eqref{eq:triangularFormDifferentLengths} the integrator chains in the $x_3$-subsystem \eqref{eq:subsys3} differ in length by one integrator, whereas in the triangular form in \cite{GstottnerKolarSchoberl:2020} those have the same length, \ie for the triangular form in \cite{GstottnerKolarSchoberl:2020}, we would have an $x_3$-subsystem of the form
		\begin{align}\label{eq:subsys3equal}
			f_3:\quad\begin{aligned}
				\dot{x}_{3,1}^1&=x_{3,1}^2&\dot{x}_{3,2}^1&=x_{3,2}^2\\
				\dot{x}_{3,1}^2&=x_{3,1}^3&\dot{x}_{3,2}^2&=x_{3,2}^3\\
				&\vdotswithin{=}&&\vdotswithin{=}\\
				\dot{x}_{3,1}^{n_3}&=u^1&\dot{x}_{3,2}^{n_3}&=u^2\,,
			\end{aligned}
		\end{align}
		instead of the form \eqref{eq:subsys3} (and $g=0$ in the $x_2$-subsystem \eqref{eq:subsys2}).
		\subsection{Characterization of the triangular form}
			In this section, we provide necessary and sufficient conditions for a two input AI-system to be static feedback equivalent to the triangular form \eqref{eq:triangularFormDifferentLengths} and thus provide a sufficient condition for such a system to be flat. Consider a two input AI-system
			\begin{align}\label{eq:aisystem}
				\begin{aligned}
					\dot{x}&=a(x)+b_1(x)u^1+b_2(x)u^2\,.
				\end{aligned}
			\end{align}
			We define the distributions $D_i$, $i=1,\ldots,n_3+1$ where $D_1=\Span{b_1,b_2}$ and $D_{i+1}=D_i+[a,D_i]$, with the smallest integer $n_3$ such that $D_{n_3+1}$ is not involutive. We again assume all distributions to have locally constant dimension and we omit discussing singularities coming along with flat outputs of \eqref{eq:triangularFormDifferentLengths} or singularities in the problem of transforming a given system into the form \eqref{eq:triangularFormDifferentLengths}. We consider generic points only, regularity conditions are omitted.
			\begin{theorem}\label{thm:1}
				The AI-system \eqref{eq:aisystem} is static feedback equivalent to the triangular form \eqref{eq:triangularFormDifferentLengths} if and only if $\Dim{D_i}=2i$, $i=1,\ldots,n_3+1$, $\C{D_{n_3+1}}\neq D_{n_3}$ and there exists a vector field $b_p=\alpha^1b_1+\alpha^2b_2$ such that with the distributions $\Delta_0=D_{n_3-1}+\Span{\Ad_a^{n_3-1}b_p}$ and $\Delta_1=D_{n_3}+\Span{\Ad_a^{n_3}b_p}$, the following conditions are satisfied:
				\begin{enumerate}[label=(\alph*)]
					\setlength{\itemsep}{5pt}
					\item $\C{\Delta_1}=\Delta_0$.\label{it:a}
					\item The derived flags of the non involutive distribution $\Delta_1$ satisfy
					\begin{align*}
						\begin{aligned}
							\Dim{\Delta_1^{(i)}}&=\Dim{\Delta_1}+i\,,\!&i&=1,\ldots,\,n_2-2,
						\end{aligned}
					\end{align*}
					with the smallest integer $n_2$ such that $\Delta_1^{(n_2-2)}=\overline{\Delta}_1$.\label{it:b}
					\item The drift satisfies the compatibility conditions\footnote{If $\overline{\Delta}_1=\mathcal{T}(\mathcal{X})$, the condition \eqref{eq:coupling} and the items \ref{it:d} and \ref{it:e} have to be omitted. In this case, the system is static feedback equivalent to \eqref{eq:triangularFormDifferentLengths} with $n_{1,1}=n_{1,2}=0$ if and only if all the other conditions of the theorem are met.}
					\begin{align}
						\begin{aligned}\label{eq:compatibilityThm1}
							[a,\C{\Delta_1^{(i)}}]&\subset\Delta_1^{(i)}\,,&i&=1,\ldots,n_2-3\,,
						\end{aligned}\\
						\begin{aligned}\label{eq:coupling}
							\Dim{\overline{\Delta}_1+[a,\Delta_1^{(n_2-3)}]}&=\Dim{\overline{\Delta}_1}+1\,.
						\end{aligned}
					\end{align}\label{it:c}
					\vspace{-1em}
					\item The distributions $G_{i+1}$, $i\geq 0$ are involutive, where $G_{i+1}=G_i+[a,G_i]$ and $G_0=\overline{\Delta}_1$.\label{it:d}
					\item $G_s=\mathcal{T}(\mathcal{X})$ holds for some integer $s$.\label{it:e}
				\end{enumerate}
			\end{theorem}\noindent
			All these conditions are easily verifiable and require differentiation and algebraic operations only. Also the construction of a vector field $b_p$, which is needed for verifying the conditions, requires differentiation and algebraic operations only, the construction is discussed in the next section. Let us outline the meaning of the individual conditions of the theorem and of the vector field $b_p$. Consider a system of the form \eqref{eq:triangularFormDifferentLengths}. The main idea of Theorem \ref{thm:1} is to characterize the three subsystems of \eqref{eq:triangularFormDifferentLengths} on their own and have separate conditions which take into account their coupling. Since the individual subsystems are either in Brunovsky normal form or essentially in extended chained form, they are in fact characterized by the Theorems \ref{thm:exactLin} and \ref{thm:extendedChainedform} in Section \ref{se:knownResults}. The $x_3$-subsystem \eqref{eq:subsys3} consists of two integrator chains with the lengths $n_3$ and $n_3-1$. The vector field $b_p$ corresponds to the input vector field of the longer integrator chain, \ie in \eqref{eq:subsys3} this would be $b_p=\partial_{x_{3,1}^{n_3}}$. The involutive distributions $D_1\subset\ldots\subset D_{n_3-1}\subset\Delta_0$ characterize the $x_3$-subsystem, \ie $\Delta_0=\Span{\partial_{x_3}}$. The top variables $x_{3,1}^1$ and $x_{3,2}^1$ of the $x_3$-subsystem act as inputs for the $x_2$-subsystem. The corresponding input vector fields of the $x_2$-subsystem read $b_{2,1}=b_1^c$ and $b_{2,2}=b_2^c+gb_1^c$ with $b_1^c=\partial_{x_2^{n_2}}$ and $b_2^c=\partial_{x_2^1}+x_2^3\partial_{x_2^2}+\ldots+x_2^{n_2}\partial_{x_2^{n_2-1}}$ being structurally of the form \eqref{eq:chaindVectorfields}. Since the shorter integrator chain of the $x_3$-subsystem has only a length of $n_3-1$, the distribution $D_{n_3}$ already contains one of these input vector field of the $x_2$-subsystem, namely the input vector field $b_{2,2}$. That is why the vector field $b_p$ is needed, it allows us to separate the $x_3$-subsystem from the $x_2$-subsystem despite the differing lengths of the integrator chains. Constructing $b_p$ roughly speaking means identifying the longer integrator chain of the $x_3$-subsystem. Given $b_p$, we can calculate the distribution $\Delta_0$ and thus explicitly identify the states which belong to the $x_3$-subsystem, \ie $\Delta_0=\Span{\partial_{x_3}}$.\\\\
			The distribution $\Delta_1$ is spanned by $\partial_{x_3}$ and both input vector fields of the $x_2$-subsystem. Item \ref{it:a} is crucial for the coupling of the $x_2$-subsystem with the $x_3$-subsystem, it assures that the $x_2$-subsystem indeed allows an AI representation with respect to its inputs $x_{3,1}^1$ and $x_{3,2}^1$. Item \ref{it:b} is in fact a condition on the distribution spanned by the input vector fields of the $x_2$-subsystem, the condition in fact matches that of Theorem \ref{thm:chainedform} for the normal chained form. Item \ref{it:b} therefore guarantees that the input vector fields of the $x_2$-subsystem can indeed be transformed into a chained structure.\\\\
			The triangular dependence of the functions $a_2^i$, $i=2,\ldots,n_2-1$, \ie the drift of the $x_2$-subsystem, on the states $x_2$ is assured by the condition \eqref{eq:compatibilityThm1}, which essentially matches the compatibility condition \eqref{eq:compatibility} in Theorem \ref{thm:extendedChainedform} for the extended chained form. The drift vector field of a system of the form \eqref{eq:triangularFormDifferentLengths} of course does not only consist of the drift vector field $a_2=a_2^2\partial_{x_2^2}+\ldots+a_2^{n_2-1}\partial_{x_2^{n_2-1}}$ of the $x_2$-subsystem, it also contains the input vector fields of the $x_2$-subsystem and has additional components which belong to the $x_1$-subsystem and the $x_3$-subsystem, the drift vector field of \eqref{eq:triangularFormDifferentLengths} actually reads $a=f_1+x_{3,2}^1(b_2^c+gb_1^c)+x_{3,1}^1b_1^c+a_2+a_3$ where $f_1=x_{1,1}^2\partial_{x_{1,1}^1}+\ldots+x_2^1\partial_{x_{1,1}^{n_{1,1}}}+x_{1,2}^2\partial_{x_{1,2}^1}+\ldots+x_2^2\partial_{x_{1,2}^{n_{1,2}}}$ and $a_3=x_{3,1}^2\partial_{x_{3,1}^1}+\ldots+x_{3,1}^{n_3}\partial_{x_{3,1}^{n_3-1}}+x_{3,2}^2\partial_{x_{3,2}^1}+\ldots+x_{3,2}^{n_3-1}\partial_{x_{3,2}^{n_3-2}}$. However, as we will see in the proof of Theorem \ref{thm:1}, the compatibility condition of the extended chained form still analogously applies, \ie when testing the compatibility of the drift of the $x_2$-subsystem, all the additional components of the drift do not matter.\\\\
			The involutive closure $\overline{\Delta}_1$ of $\Delta_1$ allows us to separate the $x_1$-subsystem from the $x_2$-subsystem and the $x_3$-subsystem, \ie $\overline{\Delta}_1=\Span{\partial_{x_3},\partial_{x_2}}$. Condition \eqref{eq:coupling} is crucial for the coupling of the $x_1$-subsystem with the $x_2$-subsystem. The static feedback linearizable $x_1$-subsystem is characterized by the involutive distributions of item \ref{it:d} and by item \ref{it:e}.\\\\
			As already mentioned, the difference between the triangular form considered here and the triangular form in \cite{GstottnerKolarSchoberl:2020} is that for the triangular form in \cite{GstottnerKolarSchoberl:2020}, we would have an $x_3$-subsystem of the form \eqref{eq:subsys3equal} instead of the form \eqref{eq:subsys3} (and $g=0$ in the $x_2$-subsystem \eqref{eq:subsys2}). An AI-system \eqref{eq:aisystem} is static feedback equivalent to this simpler triangular form presented in \cite{GstottnerKolarSchoberl:2020} if and only if the items \ref{it:a} to \ref{it:e} of Theorem \ref{thm:1} are met with $D_{n_3}$ and $D_{n_3+1}$ instead of $\Delta_0$ and $\Delta_1$. So in this case, because of the equal length of the integrator chains, the $x_3$-subsystem is simply characterized by the involutive distributions $D_1\subset\ldots\subset D_{n_3}$ and the distribution $D_{n_3+1}$ is spanned by $\partial_{x_3}$ and the input vector fields of the $x_2$-subsystem. No extra effort is needed for identifying the states which belong to the $x_3$-subsystem. Despite the similarity of the two triangular forms, the triangular form with differing lengths of the integrator chains in the $x_3$-subsystem, \ie the triangular form covered in this paper, is applicable to many practical and academic examples which cannot be handled with the triangular form presented in \cite{GstottnerKolarSchoberl:2020}. Among others, \eg the already mentioned planar VTOL aircraft and the model of a gantry crane and some academic examples considered in \cite{Schoberl:2014}, see Section \ref{se:ex}. Although the proof of Theorem \ref{thm:1} is also somewhat similar to the proof of the main theorem in \cite{GstottnerKolarSchoberl:2020}, the differing lengths of the integrator chains in the $x_3$-subsystem greatly increase the complexity of the proof. In contrast to \cite{GstottnerKolarSchoberl:2020}, where proofs are only sketched, detailed proofs are provided in this paper.
			\subsubsection{Determining a vector field $b_p$}
				According to Theorem \ref{thm:1}, the test for static feedback equivalence of an AI-system \eqref{eq:aisystem} to the triangular form \eqref{eq:triangularFormDifferentLengths} involves finding a certain linear combination $b_p=\alpha^1b_1+\alpha^2b_2$ of its input vector fields $b_1$ and $b_2$. It should be noted that in Theorem \ref{thm:1}, only the direction of $b_p$ matters, \ie if the conditions of Theorem \ref{thm:1} are met with $b_p$, then they are also met with $\tilde{b}_p=\lambda b_p$ with an arbitrary non-zero function $\lambda$ of the state of the system. A detailed proof of this property is provided in Appendix \ref{ap:detailsTheorem1}. As already mentioned, provided that the system under consideration is indeed static feedback equivalent to the triangular form \eqref{eq:triangularFormDifferentLengths}, such a vector field $b_p$ corresponds to the input vector field of the longer integrator chain in a corresponding $x_3$-subsystem. In the following, we explain the construction of such a vector field $b_p$. We will start by deriving a necessary condition, which every such vector field $b_p$ has to meet. This will allow us to determine candidates for the vector field $b_p$. We will then show, that this necessary condition yields only at most two non-collinear candidates for $b_p$. Thus, if the system is indeed static feedback equivalent to \eqref{eq:triangularFormDifferentLengths}, at least with one of those candidates, the conditions of Theorem \ref{thm:1} must indeed be met. Let us introduce the abbreviations $v_1=\Ad_a^{n_3-1}b_1$, $v_2=\Ad_a^{n_3-1}b_2$ and $v_p=\alpha^1v_1+\alpha^2v_2$. With those, because of $\Ad_a^{n_3-1}b_p=\alpha^1\Ad_a^{n_3-1}b_1+\alpha^2\Ad_a^{n_3-1}b_2\Mod D_{n_3-1}$ and analogously $\Ad_a^{n_3}b_p=\alpha^1\Ad_a^{n_3}b_1+\alpha^2\Ad_a^{n_3}b_2\Mod D_{n_3}$, we have $\Ad_a^{n_3-1}b_p=v_p\Mod D_{n_3-1}$ and $\Ad_a^{n_3}b_p=[a,v_p]\Mod D_{n_3}$. Thus, for the distributions $\Delta_0$ and $\Delta_1$ of Theorem \ref{thm:1}, we have $\Delta_0=D_{n_3-1}+\Span{v_p}$ and $\Delta_1=D_{n_3}+\Span{[a,v_p]}$. Item \ref{it:a} of Theorem \ref{thm:1} requires that $\C{\Delta_1}=\Delta_0$. Since $v_p\in\Delta_0$, this implies that $[v_p,\Delta_1]\subset\Delta_1$ must hold and in particular $[v_p,[a,v_p]]\in\Delta_1$ must hold. A necessary condition on $v_p=\alpha^1v_1+\alpha^2v_2$ (\ie a necessary condition on the coefficients $\alpha^1$ and $\alpha^2$ in this linear combination) is thus
				\begin{align}\label{eq:pde}
					\begin{aligned}
						[v_p,[a,v_p]]&=(\alpha^1)^2[v_1,[a,v_1]]+\alpha^1\Lie_{v_1}\alpha^1[a,v_1]\\
						&~~~+\alpha^1\alpha^2[v_1,[a,v_2]]+\alpha^1\Lie_{v_1}\alpha^2[a,v_2]\\
						&~~~+\alpha^2\alpha^1[v_2,[a,v_1]]+\alpha^2\Lie_{v_2}\alpha^1[a,v_1]\\
						&~~~+(\alpha^2)^2[v_2,[a,v_2]]+\alpha^2\Lie_{v_2}\alpha^2[a,v_2]\Mod D_{n_3}~\overset{!}{\in}\Delta_1\,,
					\end{aligned}
				\end{align}
				which is a system of non-linear PDEs in $\alpha^1$ and $\alpha^2$. From solutions, we obtain candidates $b_p=\alpha^1b_1+\alpha^2b_2$ (or $v_p=\alpha^1v_1+\alpha^2v_2$). Since by construction we have $\Delta_1\subset D_{n_3+1}$, a solution of \eqref{eq:pde} also meets $[v_p,[a,v_p]]\in D_{n_3+1}$. Thus, we have the weaker necessary condition $[v_p,[a,v_p]]\overset{!}{\in}D_{n_3+1}$, which because of $[a,v_1],[a,v_2]\in D_{n_3+1}$, is purely algebraic and reads\footnote{Here, we used $[v_2,[a,v_1]]=[v_1,[a,v_2]]\Mod D_{n_3+1}$, which follows from the Jacobi identity
				\begin{align*}
					\begin{aligned}
						[v_2,[a,v_1]]+\underbrace{[v_1,[v_2,a]]}_{=-[v_1,[a,v_2]]}+\underbrace{[a,\underbrace{[v_1,v_2]}_{\in D_{n_3}}]}_{\in D_{n_3+1}}&=0\,.
					\end{aligned}
				\end{align*}
				}
				\begin{align}\label{eq:alg}
					\begin{aligned}
						(\alpha^1)^2[v_1,[a,v_1]]+2\alpha^1\alpha^2[v_1,[a,v_2]]+(\alpha^2)^2[v_2,[a,v_2]]&\overset{!}{\in}D_{n_3+1}\,.
					\end{aligned}
				\end{align}
				Although \eqref{eq:alg} is again only a necessary condition which the coefficients $\alpha^1$ and $\alpha^2$ must fulfill, this necessary condition yields at most two non-collinear candidates for the vector field $b_p=\alpha^1b_1+\alpha^2b_2$. This is shown in detail in Appendix \ref{ap:analysisOfTheEquationsystem}. Thus, if the system under consideration is indeed static feedback equivalent to \eqref{eq:triangularFormDifferentLengths}, at least with one of those two non-collinear candidates for $b_p$, the conditions of Theorem \ref{thm:1} must indeed be met. 
				\begin{remark}\label{rem:simplerMethod}
					There exists a simpler method for determining a vector field $b_p$ in certain cases. If for an AI-system \eqref{eq:aisystem} $\Ad_a^{n_3+1}b_1\notin H$ or $\Ad_a^{n_3+1}b_2\notin H$, with the distribution $H=D_{n_3+1}+[D_{n_3},D_{n_3+1}]$ holds, the direction of $b_p=\alpha^1b_1+\alpha^2b_2$ is uniquely determined by the condition $\Ad_a^{n_3+1}b_p\in H$. Since $\Ad_a^{n_3+1}b_p=\alpha^1\Ad_a^{n_3+1}b_1+\alpha^2\Ad_a^{n_3+1}b_2\Mod D_{n_3+1}$ and $D_{n_3+1}\subset H$, the condition $\Ad_a^{n_3+1}b_p\in H$ yields a system of linear equations for determining $\alpha^1$ and $\alpha^2$. If $\Ad_a^{n_3+1}b_1,\,\Ad_a^{n_3+1}b_2\in H$, this method for determining $b_p$ is of course not applicable, since $\Ad_a^{n_3+1}b_p\in H$ would be met for any linear combination $b_p=\alpha^1b_1+\alpha^2b_2$ of the input vector fields of the system. A proof of this property can be found in Appendix \ref{ap:proofSimplifiedMethod}.
				\end{remark}		
		\subsection{Determining flat outputs}\label{se:determinifFlatOutputs}
			For determining flat outputs of a system which is static feedback equivalent to \eqref{eq:triangularFormDifferentLengths}, there is no need to actually transform the system into the form \eqref{eq:triangularFormDifferentLengths}. A system which is static feedback equivalent to \eqref{eq:triangularFormDifferentLengths}, meets the conditions of Theorem \ref{thm:1}. Flat outputs can be derived directly from the distributions $\Delta_1$ and $G_i$, which are involved in Theorem \ref{thm:1}, items \ref{it:b} and \ref{it:d}. All what follows in this section is actually contained in the sufficiency part of the proof of Theorem \ref{thm:1}. Here, we only summarize the computation of compatible flat outputs, for details, we refer to the proof of Theorem \ref{thm:1}. Consider a system of the form \eqref{eq:triangularFormDifferentLengths}. Depending on the length of the integrator chains in the $x_1$-subsystem \eqref{eq:subsys1}, flat outputs are determined differently. In particular, we have to distinguish between the cases that
			\begin{enumerate}[label=\arabic*)]
				\item both integrator chains have at least length one, \ie $n_{1,1},n_{1,2}\geq 1$,
				\item one of the integrator chains has length zero, and the other one at least a length of one,
				\item both integrator chains have length zero, \ie $n_{1,1}=n_{1,2}=0$, an $x_1$-subsystem does not exist at all.
			\end{enumerate}
			Given a system which meets the conditions of Theorem \ref{thm:1}, we can easily test which case applies. If we have $\Dim{G_1}=\Dim{\overline{\Delta}_1}+2$, in a corresponding triangular form \eqref{eq:triangularFormDifferentLengths}, both chains have at least length one. If $\Dim{G_1}=\Dim{\overline{\Delta}_1}+1$, one chain has length zero, if $\overline{\Delta}_1=\mathcal{T}(\mathcal{X})$, both have length zero. In the following, we discuss these three cases in more detail. 
			\paragraph*{Case 1:} If $n_{1,1},n_{1,2}\geq 1$, \ie both integrator chains of the $x_1$-subsystem \eqref{eq:subsys1} have at least length one, flat outputs are all pairs of functions $(\varphi^1,\varphi^2)$, which form a linearizing output of the $x_1$-subsystem. The $x_1$-subsystem is characterized by the distributions $G_i$ of Theorem \ref{thm:1} item \ref{it:d}. So in this case, flat outputs are determined from the sequence of involutive distributions $G_i$ of Theorem \ref{thm:1} item \ref{it:d}, in the same way as linearizing outputs are determined from the sequence of involutive distributions involved in the test for static feedback linearizability (see \eg \cite{JakubczykRespondek:1980}, \cite{NijmeijervanderSchaft:1990}).
			\paragraph*{Case 2:} Here, the $x_1$-subsystem determines one component $\varphi^1$ of a flat output. This function is obtained by integrating $G_{s-1}^\perp$, \ie by finding a function $\varphi^1$ such that $\Span{\D\varphi^1}=G_{s-1}^\perp$, with $G_{s-1}$ of Theorem \ref{thm:1}, item \ref{it:d}. A possible second component $\varphi^2$ is obtained by integrating the integrable codistribution $L^\perp=(\Delta_1^{(n_2-3)})^\perp+\Span{\D\Lie_a^s\varphi^1}$. A function $\varphi^2$ whose differential $\D\varphi^2$ together with $\{\D\varphi^1,\D\Lie_a\varphi^1,\ldots,\D\Lie_a^s\varphi^1\}$ spanns the codistribution $L^\perp$, is a possible second component.
			\paragraph*{Case 3:} If both chains have length zero (the $x_1$-subsystem does not exist at all), the problem of finding flat outputs, is in fact the same as finding flat outputs of a system that is static feedback equivalent to the chained form. This problem is addressed in \cite{LiRespondek:2012}. In this case, flat outputs are all pairs of functions $(\varphi^1,\varphi^2)$, which meet $L=(\Span{\D\varphi^1,\D\varphi^2})^\perp\subset\Delta_1^{(n_2-3)}$, with $\Delta_1^{(n_2-3)}$ of Theorem \ref{thm:1} item \ref{it:b}. In \cite{LiRespondek:2012}, Theorem 2.10, a method for constructing such a distribution $L$ is provided. The distribution $L$ is not unique, one has to choose one function $\varphi^1$ whose differential $\D\varphi^1\neq 0$ annihilates $\mathcal{C}(\Delta_1^{(n_2-3)})$. Once such a function has been chosen, the distribution $L$ can be calculated, and in turn a possible second function $\varphi^2$, which together with $\varphi^1$ forms a possible flat output, can be calculated. Equivalent to the method for determining $L$ provided in \cite{LiRespondek:2012}, once a function $\varphi^1$ whose differential annihilates $\mathcal{C}(\Delta_1^{(n_2-3)})$ has been chosen, the annihilator of $L$ can also be calculated via $L^\perp=(\Delta_1^{(n_2-3)})^\perp+\Span{\D\varphi^1}$.\\\\
			Note that Case 2 and 3 are in fact similar. The function $\Lie_a^s\varphi^1$ in Case 2 corresponds to $\varphi^1$ from Case 3. In Case 3, we have to choose the function $\varphi^1$. Once this function has been chosen, the distribution $L$ is uniquely determined by this function. In Case 2, the distribution $L$ is uniquely determined by the function $\Lie_a^s\varphi^1$, which is imposed by the $x_1$-subsystem.
		\subsection{Proof of Theorem \ref{thm:1}}
			To keep the proof reasonably compact, parts of it are condensed into propositions and small facts, which are proven in Appendix \ref{ap:detailsTheorem1}. The following two lemmas are of particular importance for the sufficiency part of the proof. Proofs of these lemmas are provided in Appendix \ref{ap:additionalProofs}.
			\begin{lemma}\label{lem:characteristic}
				Let $D$ be a distribution. Every characteristic vector field of $D$, \ie every vector field $c\in\C{D}$ is also characteristic for its derived flag $D^{(1)}$, \ie $\C{D}\subset\C{D^{(1)}}$.
			\end{lemma}\noindent
			An immediate consequence of Lemma \ref{lem:characteristic} is that the Cauchy characteristic distributions $\C{D^{(i)}}$, $i\geq 0$ form the sequence of nested involutive distributions $\C{D}\subset\C{D^{(1)}}\subset\C{D^{(2)}}\subset\ldots$
			\begin{lemma}\label{lem:cartan}
				If a $d$-dimensional distribution $D$ satisfies $\Dim{\C{D}}=d-2$ and $\Dim{D^{(i)}}=d+i$, $i=1,\ldots,l$ with $l$ such that $D^{(l)}=\overline{D}$, then the Cauchy characteristics $\C{D^{(i)}}$ satisfy $\Dim{\C{D^{(i)}}}=d-2+i$ and $\C{D^{(i)}}\subset D^{(i-1)}$, $i=1,\ldots,l-1$.
			\end{lemma}\noindent
			Lemma \ref{lem:cartan} is based on a similar one in \cite{Cartan:1914}, see also \cite{MartinRouchon:1994}, Lemma 2.\\\\
			\textit{\textbf{Necessity.}} We have to show that a system of the form \eqref{eq:triangularFormDifferentLengths} meets the conditions of Theorem \ref{thm:1}. Recall that the drift vector field of a system of the form \eqref{eq:triangularFormDifferentLengths} reads $a=f_1+x_{3,2}^1(b_2^c+gb_1^c)+x_{3,1}^1b_1^c+a_2+a_3$ where $b_1^c=\partial_{x_2^{n_2}}$, $b_2^c=\partial_{x_2^1}+x_2^3\partial_{x_2^2}+\ldots+x_2^{n_2}\partial_{x_2^{n_2-1}}$, $f_1=x_{1,1}^2\partial_{x_{1,1}^1}+\ldots+x_2^1\partial_{x_{1,1}^{n_{1,1}}}+x_{1,2}^2\partial_{x_{1,2}^1}+\ldots+x_2^2\partial_{x_{1,2}^{n_{1,2}}}$ and $a_3=x_{3,1}^2\partial_{x_{3,1}^1}+\ldots+x_{3,1}^{n_3}\partial_{x_{3,1}^{n_3-1}}+x_{3,2}^2\partial_{x_{3,2}^1}+\ldots+x_{3,2}^{n_3-1}\partial_{x_{3,2}^{n_3-2}}$. The input vector fields of \eqref{eq:triangularFormDifferentLengths} are given by $b_1=\partial_{x_{3,1}^{n_3}}$ and $b_2=\partial_{x_{3,2}^{n_3-1}}$. The distributions defined right before Theorem \ref{thm:1} are thus given by
			\begin{align*}
				\begin{aligned}
					D_1&=\Span{\partial_{x_{3,1}^{n_3}},\partial_{x_{3,2}^{n_3-1}}}\\
					D_2&=\Span{\partial_{x_{3,1}^{n_3}},\partial_{x_{3,2}^{n_3-1}},\partial_{x_{3,1}^{n_3-1}},\partial_{x_{3,2}^{n_3-2}}}\\
					&\vdotswithin{=}\\
					D_{n_3-1}&=\Span{\partial_{x_{3,1}^{n_3}},\partial_{x_{3,2}^{n_3-1}},\ldots,\partial_{x_{3,1}^2},\partial_{x_{3,2}^1}}\\
					D_{n_3}&=\Span{\underbrace{\partial_{x_{3,1}^{n_3}},\partial_{x_{3,2}^{n_3-1}},\ldots,\partial_{x_{3,1}^2},\partial_{x_{3,2}^1},\partial_{x_{3,1}^1}}_{\partial_{x_3}},b_2^c+gb_1^c}\\
					D_{n_3+1}&=\Span{\partial_{x_3},b_1^c,b_2^c,[a,b_2^c+gb_1^c]}\,.
				\end{aligned}
			\end{align*}
			The distributions $D_1,\ldots,D_{n_3}$ are involutive and they meet $\Dim{D_i}=2i$. 
			\begin{fact}\label{ft:1}
				The distribution $D_{n_3+1}$ is not involutive, it meets $\Dim{D_{n_3+1}}=2n_3+2$ and the condition $\C{D_{n_3+1}}\neq D_{n_3}$ is met. 
			\end{fact}
			The input vector field belonging to the longer integrator chain in the $x_3$-subsystem is $b_1=\partial_{x_{3,1}^{n_3}}$. With $b_p=b_1$, we obtain $\Ad_a^{n_3-1}b_p=(-1)^{n_3-1}\,\partial_{x_{3,1}^1}$ and $\Ad_a^{n_3}b_p=(-1)^{n_3}\,b_1^c$, and thus $\Delta_0=D_{n_3-1}+\Span{\Ad_a^{n_3-1}b_p}=\Span{\partial_{x_3}}$ and  $\Delta_1=D_{n_3}+\Span{\Ad_a^{n_3}b_p}=\Span{\partial_{x_3},b_1^c,b_2^c}$. With these two distributions, the items \ref{it:a} to \ref{it:e} of Theorem \ref{thm:1} are met. Item \ref{it:a} is met since the vector fields $b_1^c$ and $b_2^c$ are independent of the variables $x_3$. We have
			\begin{align*}
				\begin{aligned}
					\Delta_1^{(i)}&=\Span{\partial_{x_3},\partial_{x_2^1}+x_2^3\partial_{x_2^2}+\ldots+x_2^{n_2-i}\partial_{x_2^{n_2-i-1}},\partial_{x_2^{n_2}},\ldots,\partial_{x_2^{n_2-i}}}\,,&i&=0,\ldots,n_2-2\,,
				\end{aligned}
			\end{align*}
			where the $(n_2-2)$-th derived flag of $\Delta_1$ is actually the involutive closure of $\Delta_1$, \ie $\Delta_1^{(n_2-2)}=\overline{\Delta}_1=\Span{\partial_{x_3},\partial_{x_2}}$, the last non involutive distribution of this sequence is $\Delta_1^{(n_2-3)}=\Span{\partial_{x_3},\partial_{x_2^1}+x_2^3\partial_{x_2^2},\partial_{x_2^{n_2}},\ldots,\partial_{x_2^3}}$. Their Cauchy characteristic distributions are given by
			\begin{align*}
				\begin{aligned}
					\C{\Delta_1^{(i)}}&=\Span{\partial_{x_3},\partial_{x_2^{n_2}},\ldots,\partial_{x_2^{n_2-i+1}}}\,,&i&=1,\ldots,n_2-3\,.
				\end{aligned}
			\end{align*}
			The distributions $\Delta_1^{(i)}$ meet the condition $\Dim{\Delta_1^{(i)}}=\Dim{\Delta_1}+i$, $i=1,\ldots,n_2-2$ of item \ref{it:b}. To show that the condition \eqref{eq:compatibilityThm1} of item \ref{it:c}, \ie $[a,\C{\Delta_1^{(i)}}]\subset\Delta_1^{(i)}$, $i=1,\ldots,n_2-3$ is met, note that we have\footnote{Evaluated for \eg $i=1$, we obtain
			\begin{align*}
				\begin{aligned}
					[\partial_{x_2^{n_2}},a]&=(x_{3,2}^1+\partial_{x_2^{n_2}}a^{n_2-1})\partial_{x_2^{n_2-1}}\Mod\C{\Delta_1^{(1)}}\,.
				\end{aligned}
			\end{align*}
			Because of $\Delta_1^{(1)}=\Span{\partial_{x_3},\partial_{x_2^1}+x_2^3\partial_{x_2^2}+\ldots+x_2^{n_2-1}\partial_{x_2^{n_2-2}},\partial_{x_2^{n_2}},\partial_{x_2^{n_2-1}}}$ and $\C{\Delta_1^{(1)}}=\Span{\partial_{x_3},\partial_{x_2^{n_2}}}$, it indeed follows that $[a,\C{\Delta_1^{(1)}}]\subset\Delta_1^{(1)}$.}
			\begin{align*}
				\begin{aligned}
					[\partial_{x_2^{n_2-i+1}},a]&=(x_{3,2}^1+\partial_{x_2^{n_2-i+1}}a^{n_2-i})\partial_{x_2^{n_2-i}}\Mod\C{\Delta_1^{(i)}}\,,&i&=1,\ldots,n_2-3\,,
				\end{aligned}
			\end{align*}
			\ie $[\partial_{x_2^{n_2-i+1}},a]\in\Delta_1^{(i)}$, $i=1,\ldots,n_2-3$. Condition \eqref{eq:coupling} of item \ref{it:c} is met since
			\begin{align*}
				\begin{aligned}
					\overline{\Delta}_1+[a,\Delta_1^{(n_2-3)}]&=\overline{\Delta}_1+\Span{\underbrace{[a,\partial_{x_2^{n_2}}],\ldots,[a,\partial_{x_2^3}]}_{\in\overline{\Delta}_1},\underbrace{[a,\partial_{x_2^1}+x_2^3\partial_{x_2^2}]}_{\notin\overline{\Delta}_1}}\,,
				\end{aligned}
			\end{align*}
			\ie $[a,\Delta_1^{(n_2-3)}]$ yields only one new direction with respect to $\overline{\Delta}_1$. For the distributions $G_i$ of item \ref{it:d}, we obtain
			\begin{align*}
				\begin{aligned}
					G_0&=\overline{\Delta}_1=\Span{\partial_{x_3},\partial_{x_2}}\\
					G_i&=\Span{\partial_{x_3},\partial_{x_2},\partial_{x_{1,1}^{n_{1,1}}},\ldots,\partial_{x_{1,1}^{n_{1,1}-i+1}},\partial_{x_{1,2}^{n_{1,2}}},\ldots,\partial_{x_{1,2}^{n_{1,2}-i+1}}}\,,&i&\geq 1\,,
				\end{aligned}
			\end{align*}
			where any $\partial_{x_{1,j}^k}$, $j=1,2$ with $k\leq 0$ has to be omitted. These distributions are obviously involutive. Thus, item \ref{it:d} is met. Item \ref{it:e} is met since we have $G_s=\mathcal{T}(\mathcal{X})$ for $s=\max\{n_{1,1},n_{1,2}\}$.\\\\
			\textit{\textbf{Sufficiency.}} We have to show that an AI-system which meets the conditions of Theorem \ref{thm:1} can be transformed into the triangular form \eqref{eq:triangularFormDifferentLengths}. Item \ref{it:a} implies that $\Delta_0$ is involutive. Because of Lemma \ref{lem:characteristic}, the Cauchy characteristics of the derived flags of $\Delta_1$ form the sequence of nested involutive distributions
			\begin{align*}
				\begin{aligned}
					\C{\Delta_1}\subset\C{\Delta_1^{(1)}}\subset\ldots\subset\C{\Delta_1^{(n_2-3)}}\,,
				\end{aligned}
			\end{align*}
			where $\C{\Delta_1}=\Delta_0$. We thus have the sequence of nested involutive distributions
			\begin{align}\label{eq:involutiveSequenceProof}
				\begin{aligned}
					&D_1\subset\ldots\subset D_{n_3-1}\subset\Delta_0\subset\C{\Delta_1^{(1)}}\subset\ldots\subset\\
					&~~~~~~~~~~\C{\Delta_1^{(n_2-3)}}\subset
					\overline{\Delta}_1\subset G_1\subset\ldots\subset G_s=\mathcal{T}(\mathcal{X})\,.
				\end{aligned}
			\end{align}
			The transformation of \eqref{eq:aisystem} into the form \eqref{eq:triangularFormDifferentLengths} is done in the following six steps.
			\paragraph*{Step 1:} Straighten out all the distributions \eqref{eq:involutiveSequenceProof} simultaneously.
			\begin{prop}\label{prop:1}
				In coordinates in which the distributions \eqref{eq:involutiveSequenceProof} are straightened out, the system \eqref{eq:aisystem} takes the form\vspace{-1em}
				\begin{align}\label{eq:step1}
					\begin{aligned}
						\dot{x}_1&=f_1(x_1,x_2^1,x_2^2,x_2^3)\\
						\dot{x}_2&=f_2(x_1,x_2,x_3^1,x_3^2,x_3^3)\\
						\dot{x}_3&=f_3(x_1,x_2,x_3,u^1,u^2)\,.
					\end{aligned}
				\end{align}
			\end{prop}
			\vspace{-3em}
			\begin{prop}\label{prop:2}
				The subsystems in \eqref{eq:step1} meet the rank conditions $\Rank{\partial_{(x_2^1,x_2^2,x_2^3)}f_1}\leq2$, $\Rank{\partial_{(x_3^1,x_3^2,x_3^3)}f_2}=2$ and $\Rank{\partial_{(x_3^2,x_3^3)}f_2}=1$\footnote{For $n_3=1$, we have $\dot{x}_2=f_2(x_1,x_2,x_3^1,u^1,u^2)$ and the latter two rank conditions read $\Rank{\partial_{(u^1,u^2,x_3^1)}f_2}=2$ and $\Rank{\partial_{(u^1,u^2)}f_2}=1$.}.
			\end{prop}
			Note that by only straightening out the distributions $\Delta_0$ and $\overline{\Delta}_1$, \ie applying a change of coordinates such that $\Delta_0=\Span{\partial_{x_3}}$ and $\overline{\Delta}_1=\Span{\partial_{x_3},\partial_{x_2}}$, we already obtain a decomposition of the system \eqref{eq:aisystem} into three subsystems. Simultaneously straightening out the remaining distributions of \eqref{eq:involutiveSequenceProof} only affects the structure of these three subsystems. In particular, by straightening out $D_1\subset\ldots\subset D_{n_3-1}$ and $G_1\subset\ldots\subset G_{s-1}$, the $x_1$-subsystem and the $x_3$-subsystem take a triangular structure, known from the static feedback linearization problem (see \eg \cite{NijmeijervanderSchaft:1990}). The inputs $u^1$ and $u^2$ of course occur affine in $f_3$, since we started with an AI-system and only applied a state transformation, which of course preserves the AI structure.
			\paragraph*{Step 2:} Transform the subsystem $\dot{x}_1=f_1(x_1,x_2^1,x_2^2,x_2^3)$ into Brunovsky normal form, \ie separate it into two integrator chains, by successively introducing new coordinates from top to bottom. In the prior to last step, we then have
			\begin{align}\label{eq:step2}
				\begin{aligned}
					f_1:\quad&\begin{aligned}
						\dot{x}_{1,1}^1&=x_{1,1}^2&\dot{x}_{1,2}^1&=x_{1,2}^2\\
						\dot{x}_{1,1}^2&=x_{1,1}^3&\dot{x}_{1,2}^2&=x_{1,2}^3\\
						&\vdotswithin{=}&&\vdotswithin{=}\\
						\dot{x}_{1,1}^{n_{1,1}}&=\varphi^1(\bar{x}_1,x_2^1,x_2^2,x_2^3)&\dot{x}_{1,2}^{n_{1,2}}&=\varphi^2(\bar{x}_1,x_2^1,x_2^2,x_2^3)
					\end{aligned}\\[1ex]
					&\quad\begin{aligned}
						\dot{x}_2&=\bar{f}_2(\bar{x}_1,x_2,x_3^1,x_3^2,x_3^3)\\
					\dot{x}_3&=\bar{f}_3(\bar{x}_1,x_2,x_3,u^1,u^2)\,,
					\end{aligned}
				\end{aligned}
			\end{align}
			with $\bar{x}_1=(x_{1,1}^1,\ldots,x_{1,1}^{n_{1,1}},x_{1,2}^1,\ldots,x_{1,2}^{n_{1,2}})$. In the following, we have to distinguish between the three possible cases regarding the actual number of integrator chains in the $x_1$-subsystem, which were already mentioned in Section \ref{se:determinifFlatOutputs}. The rank of the Jacobian matrix $\partial_{(x_2^1,x_2^2,x_2^3)}f_1$, corresponds to the actual number of non-redundant inputs of the $x_1$-subsystem and thus to the actual number of integrator chains in the $x_1$-subsystem, which can either be two, one or zero.
			\subparagraph*{Case 1:} If $\Rank{\partial_{(x_2^1,x_2^2,x_2^3)}f_1}=2$ holds, the functions $\varphi^j(\bar{x}_1,x_2^1,x_2^2,x_2^3)$, $j=1,2$ in \eqref{eq:step2} determine the desired top variables for the $x_2$-subsystem. Because of $\Rank{\partial_{(x_2^1,x_2^2,x_2^3)}f_1}=2$, these functions meet $\D \bar{x}_1\wedge\D\varphi^1\wedge\D\varphi^2\neq 0$ and thus, they can indeed serve as states for the $x_2$-subsystem. For the system to be static feedback equivalent to the triangular form \eqref{eq:triangularFormDifferentLengths}, these functions have to form a flat output of the $x_2$-subsystem, which is compatible with its (extended) chained structure. For that, these functions have to meet $L=(\Span{\D\bar{x}_1,\D\varphi^1,\D\varphi^2})^\perp\subset\Delta_1^{(n_2-3)}$, which in turn implies that they indeed form a flat output compatible with the (extended) chained form of the $x_2$-subsystem (the distribution $L$ is of importance in the problem of transforming the $x_2$-subsystem into (extended) chained form). 
			\begin{prop}\label{prop:3}
				If $\Rank{\partial_{(x_2^1,x_2^2,x_2^3)}f_1}=2$ holds, the functions $\varphi^j(\bar{x}_1,x_2^1,x_2^2,x_2^3)$, $j=1,2$ in \eqref{eq:step2} satisfy $L=(\Span{\D\bar{x}_1,\D\varphi^1,\D\varphi^2})^\perp\subset\Delta_1^{(n_2-3)}$.
			\end{prop}
			\vspace{-1.7em}
			\subparagraph*{Case 2:} If $\Rank{\partial_{(x_2^1,x_2^2,x_2^3)}f_1}=1$ holds, the $x_1$-subsystem only consists of one integrator chain (one chain in \eqref{eq:step2} is missing), so it determines only one function $\varphi^1(\bar{x}_1,x_2^1,x_2^2,x_2^3)$ which we want as top variable in the $x_2$-subsystem. In this case, there always exists a second function $\varphi^2(\bar{x}_1,x_2^1,x_2^2,x_2^3)$ which together with $\varphi^1$ fulfills $L=(\Span{\D \bar{x}_1,\D\varphi^1,\D\varphi^2})^\perp\subset\Delta_1^{(n_2-3)}$.
			\begin{prop}\label{prop:4}
				If $\Rank{\partial_{(x_2^1,x_2^2,x_2^3)}f_1}=1$ holds, there always exists a function $\varphi^2(\bar{x}_1,x_2^1,x_2^2,x_2^3)$ which together with $\varphi^1$ fulfills $L=(\Span{\D \bar{x}_1,\D\varphi^1,\D\varphi^2})^\perp\subset\Delta_1^{(n_2-3)}$.
			\end{prop}
			\vspace{-1.7em}
			\subparagraph*{Case 3:} Finally, in the case that the $x_1$-subsystem does not exist at all, \ie $\overline{\Delta}_1=\mathcal{T}(\mathcal{X})$, two functions $\varphi^j(x_2^1,x_2^2,x_2^3)$, which fulfill $L=(\Span{\D\varphi^1,\D\varphi^2})^\perp\subset\Delta_1^{(n_2-3)}$, have to be found. Such functions always exist.
			\begin{prop}\label{prop:5}
				If $\overline{\Delta}_1=\mathcal{T}(\mathcal{X})$, there always exist two functions $\varphi^j(x_2^1,x_2^2,x_2^3)$, which fulfill $L=(\Span{\D\varphi^1,\D\varphi^2})^\perp\subset\Delta_1^{(n_2-3)}$.
			\end{prop}
			\vspace{-1.7em}
			\paragraph*{Step 3:} Introduce the functions $\varphi^j(\bar{x}_1,x_2^1,x_2^2,x_2^3)$, as the top variables of the $x_2$-subsystem, \ie apply the state transformation $\tilde{x}_2^j=\varphi^j(\bar{x}_1,x_2^1,x_2^2,x_2^3)$, $j=1,2$. (When introducing these function as new states of the $x_2$-subsystem, it may be necessary to also introduce $\tilde{x}_2^3=\varphi^3(\bar{x}_1,x_2^1,x_2^2,x_2^3)$ with a suitably chosen function $\varphi^3$, since it may happen that $\D \bar{x}_1\wedge\D\varphi^1\wedge\D\varphi^2\wedge\D x_2^3=0$.) This completes the transformation of the $x_1$-subsystem into Brunovsky normal form, \ie
			\begin{align*}
				\begin{aligned}
					f_1:\quad&\begin{aligned}
					\dot{x}_{1,1}^1&=x_{1,1}^2&\dot{x}_{1,2}^1&=x_{1,2}^2\\
					&\vdotswithin{=}&&\vdotswithin{=}\\
					\dot{x}_{1,1}^{n_{1,1}}&=\tilde{x}_2^1&\dot{x}_{1,2}^{n_{1,2}}&=\tilde{x}_2^2
					\end{aligned}\\[1ex]
					&\quad\begin{aligned}
					\dot{\tilde{x}}_2&= \tilde{f}_2(\bar{x}_1,\tilde{x}_2,x_3^1,x_3^2,x_3^3)\\
					\dot{x}_3&=\tilde{f}_3(\bar{x}_1,\tilde{x}_2,x_3,u^1,u^2)\,,
					\end{aligned}
				\end{aligned}
			\end{align*}
			where $\tilde{x}_2=(\tilde{x}_2^1,\tilde{x}_2^2,x_2^3,\ldots,x_2^{n_2})$. Furthermore, this transformation straightens out the distribution $L=(\Span{\D \bar{x}_1,\D\varphi^1,\D\varphi^2})^\perp$ simultaneously with the distributions \eqref{eq:involutiveSequenceProof}, \ie $L=\Span{\partial_{x_3},\partial_{x_2^{n_2}},\ldots,\partial_{x_2^3}}$.\\\\
			The following two steps deal with the transformation of the $x_2$-subsystem into essentially (extended) chained form, see also Remark \ref{rem:essentiallyExtendedChained}.
			\paragraph*{Step 4:} Because of the rank condition $\Rank{\partial_{(x_3^2,x_3^3)}f_2}=1$ of Step 1, we also have $\Rank{\partial_{(x_3^2,x_3^3)}\tilde{f}_2}=1$. Therefore, without loss of generality, we can assume that $\tilde{f}_2$ explicitly depends on $x_3^2$, if not, swap $x_3^2$ and $x_3^3$. Let us assume that the first component $\tilde{f}_2^1$ of $\tilde{f}_2$ explicitly depends on $x_3^2$ (after eventually swapping $x_3^2$ and $x_3^3$). This enables us to replace the state $x_3^2$ of the $x_3$-subsystem by the new state\footnote{In case that $n_3=1$ holds, instead of the states $x_3^2$ or $x_3^3$ at least one of the inputs $u^1$ or $u^2$ occurs in $\tilde{f}_2$ (the inputs $u^1$ and $u^2$ would of course occur affine in $\tilde{f}_2$). Instead of the state transformation \eqref{eq:firstInput}, we then have the input transformation $\tilde{u}^2=\tilde{f}_2^1(\bar{x}_1,\tilde{x}_2,x_3^1,u^1,u^2)$ and $u^1$ left unchanged. Crucial for this transformation to be regular is that $\tilde{f}_2^1$ indeed depends on $u^2$ (after eventually swapping $u^1$ and $u^2$).\label{fn:transformation}}
			\begin{align}\label{eq:firstInput}
				\begin{aligned}
					x_{3,2}^1&=\tilde{f}_2^1(\bar{x}_1,\tilde{x}_2,x_3^1,x_3^2,x_3^3)\,,
				\end{aligned}
			\end{align}
			\ie with this new state, we replace $x_3^2$ and leave all the other coordinates unchanged. This transformation normalizes the first equation of the $x_2$-subsystem, \ie $\dot{\tilde{x}}_2^1=x_{3,2}^1$. The following fact guarantees that this transformation is indeed a regular transformation.
			\begin{fact}\label{ft:2}
				The function $\tilde{f}_2^1$ in \eqref{eq:firstInput} indeed explicitly depends on $x_3^2$ (after eventually swapping $x_3^2$ and $x_3^3$).
			\end{fact}
			\vspace{-2em}
			\begin{prop}\label{prop:6}
				After applying the transformation \eqref{eq:firstInput}, the $x_2$-subsystem reads
				\begin{align}\label{eq:subsys2triangular}
					f_2:\quad\begin{aligned}
						\dot{\tilde{x}}_2^1&=x_{3,2}^1\\
						\dot{\tilde{x}}_2^2&=b_2^2(\bar{x}_1,\tilde{x}_2^1,\tilde{x}_2^2,x_2^3)x_{3,2}^1+a_2^2(\bar{x}_1,\tilde{x}_2^1,\tilde{x}_2^2,x_2^3)\\
						\dot{x}_2^3&=b_2^3(\bar{x}_1,\tilde{x}_2^1,\tilde{x}_2^2,x_2^3,x_2^4)x_{3,2}^1+a_2^3(\bar{x}_1,\tilde{x}_2^1,\tilde{x}_2^2,x_2^3,x_2^4)\\
						&\vdotswithin{=}\\
						\dot{x}_2^{n_2-1}&=b_2^{n_2-1}(\bar{x}_1,\tilde{x}_2)x_{3,2}^1+a_2^{n_2-1}(\bar{x}_1,\tilde{x}_2)\\
						\dot{x}_2^{n_2}&=g(\bar{x}_1,\tilde{x}_2,x_3^1,x_{3,2}^1)\,.
					\end{aligned}
				\end{align}
			\end{prop}
			The triangular dependence of the functions $b_2^i$, $i=2,\ldots,n_2-1$ on the states $(x_2^3,\ldots,x_2^{n_2})$ is in fact a consequence of item \ref{it:b}. The triangular dependence of the functions $a_2^i$, $i=2,\ldots,n_2-1$ on these states is guaranteed by item \ref{it:c} condition \eqref{eq:compatibilityThm1}, \ie $[a,\C{\Delta_1^{(i)}}]\subset\Delta_1^{(i)}$, $i=1,\ldots,n_2-3$. 
			\paragraph*{Step 5:} Successively introduce the functions $b_2^i$ as new states in the $x_2$-subsystem from top to bottom. After $n_2-2$ such steps, the $x_2$-subsystem reads
			\begin{align}\label{eq:step5}
				f_2:\quad\begin{aligned}
					\dot{\tilde{x}}_2^1&=x_{3,2}^1\\
					\dot{\tilde{x}}_2^2&=\tilde{x}_2^3x_{3,2}^1+\tilde{a}_2^2(\bar{x}_1,\tilde{x}_2^1,\tilde{x}_2^2,\tilde{x}_2^3)\\
					\dot{\tilde{x}}_2^3&=\tilde{x}_2^4x_{3,2}^1+\tilde{a}_2^3(\bar{x}_1,\tilde{x}_2^1,\tilde{x}_2^2,\tilde{x}_2^3,\tilde{x}_2^4)\\
					&\vdotswithin{=}\\
					\dot{\tilde{x}}_2^{n_2-1}&=\tilde{x}_2^{n_2}x_{3,2}^1+\tilde{a}_2^{n_2-1}(\bar{x}_1,\tilde{x}_2)\\
					\dot{\tilde{x}}_2^{n_2}&=\tilde{g}(\bar{x}_1,\tilde{x}_2,x_3^1,x_{3,2}^1)\,.
				\end{aligned}
			\end{align}
			\begin{remark}\label{rem:complete}
				Introducing $x_{3,1}^1=\tilde{g}(\bar{x}_1,\tilde{x}_2,x_3^1,x_{3,2}^1)$ would complete the transformation of the $x_2$-subsystem to extended chained form (except for the dependence of the drift $\tilde{a}_2$ on the states $\bar{x}_1$, see also Remark \ref{rem:essentiallyExtendedChained}). However, after this transformation, the distribution $D_{n_3-1}$ would in general no longer be straightened out.
			\end{remark}\noindent
			From the involutivity of $D_{n_3}$, it follows that in the last line of \eqref{eq:step5}, we actually have $\dot{\tilde{x}}_2^{n_2}=\tilde{g}^1(\bar{x}_1,\tilde{x}_2,x_3^1)+\tilde{g}^2(\bar{x}_1,\tilde{x}_2)x_{3,2}^1$, \ie $\partial_{x_{3,2}^1}^2\tilde{g}=0$ and $\partial_{x_3^1}\partial_{x_{3,2}^1}\tilde{g}=0$\footnote{Calculating $D_{n_3}$ in the coordinates obtained so far, we obtain $D_{n_3}=\Span{\partial_{x_3},\partial_{\tilde{x}_2^1}+\tilde{x}_2^3\partial_{\tilde{x}_2^2}+\ldots+\tilde{x}_2^{n_2}\partial_{\tilde{x}_2^{n_2-1}}+\partial_{x_{3,2}^1}\tilde{g}\partial_{\tilde{x}_2^{n_2}}}$.  If the function $\tilde{g}$ would be of the general nonlinear form $\tilde{g}(\bar{x}_1,\tilde{x}_2,x_3^1,x_{3,2}^1)$ with $\partial_{x_{3,2}^1}^2\tilde{g}\neq 0$ or $\partial_{x_3^1}\partial_{x_{3,2}^1}\tilde{g}\neq 0$, the distribution $D_{n_3}$ would not be involutive, which contradicts with $D_{n_3}$ being indeed involutive.}. Introducing $x_{3,1}^1=\tilde{g}^1(\bar{x}_1,\tilde{x}_2,x_3^1)$, \ie replacing $x_3^1$ by the new state $x_{3,1}^1$ and leaving all the other coordinates unchanged, keeps all the distributions \eqref{eq:involutiveSequenceProof} straightened out and results in $\dot{\tilde{x}}_2^{n_2}=x_{3,1}^1+\tilde{g}^2(\bar{x}_1,\tilde{x}_2)x_{3,2}^1$.
			\paragraph*{Step 6:} Transform the $x_3$-subsystem into Brunovsky normal form, by successively introducing new states from top to bottom and finally applying a suitable static feedback.\\\\
			The transformation of an AI-system into the triangular form \eqref{eq:triangularFormDifferentLengths} by following these six steps is demonstrated on two different examples in the following section.
	\section{Examples}\label{se:ex}
		\subsection{Planar VTOL aircraft}
			Consider again our motivating example, the planar VTOL aircraft \eqref{eq:vtol}. In the following, we apply Theorem \ref{thm:1} to show that this system is indeed static feedback equivalent to the triangular form \eqref{eq:triangularFormDifferentLengths}. Based on that, we will derive a possible flat output compatible with the triangular form. Finally, we will explicitly transform the system into the form \eqref{eq:triangularFormDifferentLengths}. The input vector fields of \eqref{eq:vtol} are given by $b_1=-\sin(\theta)\partial_{v_x}+\cos(\theta)\partial_{v_z}$ and $b_2=\epsilon\cos(\theta)\partial_{v_x}+\epsilon\sin(\theta)\partial_{v_z}+\partial_\omega$. The drift is given by $a=v_x\partial_x+v_z\partial_z+\omega\partial_\theta-\partial_{v_z}$. The distribution $D_1=\Span{b_1,b_2}$ is involutive, the distribution
			\begin{align*}
				D_2&=D_1+[a,D_1]\\
				&=\Span{-\sin(\theta)\partial_{v_x}+\cos(\theta)\partial_{v_z},\epsilon\cos(\theta)\partial_{v_x}+\epsilon\sin(\theta)\partial_{v_z}+\partial_\omega,\sin(\theta)\partial_x-\cos(\theta)\partial_z-\\
				&\hspace{3em}\omega\cos(\theta)\partial_{v_x}-\omega\sin(\theta)\partial_{v_z},\epsilon\cos(\theta)\partial_x+\epsilon\sin(\theta)\partial_z+\partial_\theta+\epsilon\omega\sin(\theta)\partial_{v_x}-\epsilon\omega\cos(\theta)\partial_{v_z}}
			\end{align*}
			is not involutive, so we have $n_3=1$ and the conditions $\Dim{D_i}=2i$, $i=1,2$ hold. The condition $\C{D_2}\neq D_1$ of Theorem \ref{thm:1} is also met. Before we can evaluate the remaining conditions of Theorem \ref{thm:1}, we have to construct a vector field $b_p$ for this system. For the distribution $H=D_2+[D_1,D_2]$, see Remark \ref{rem:simplerMethod}, we obtain
			\begin{align*}
				\begin{aligned}
					H&=\Span{\sin(\theta)\partial_x-\cos(\theta)\partial_z,\epsilon\partial_x+\cos(\theta)\partial_\theta,\partial_{v_x},\partial_{v_z},\partial_\omega}
				\end{aligned}
			\end{align*}
			and for the vector fields $\Ad_a^{n_3+1}b_1$, $\Ad_a^{n_3+1}b_2$ we have
			\begin{align*}
				\begin{aligned}
					\Ad_a^2b_1&=2\omega\cos(\theta)\partial_x+2\omega\sin(\theta)\partial_z+\omega^2\sin(\theta)\partial_{v_x}-\omega^2\cos(\theta)\partial_{v_z}\,,\\
					\Ad_a^2b_2&=2\epsilon\omega\sin(\theta)\partial_x-2\epsilon\omega\cos(\theta)\partial_z-\epsilon\omega^2\cos(\theta)\partial_{v_x}-\epsilon\omega^2\sin(\theta)\partial_{v_z}\,.
				\end{aligned}
			\end{align*}
			The vector field $\Ad_a^2b_2$ is contained in $H$, the vector field $\Ad_a^2b_1$ is not contained in $H$, so we have $b_p=b_2$ (see again Remark \ref{rem:simplerMethod}). Thus, for the distributions $\Delta_0=D_{n_3-1}+\Span{\Ad_a^{n_3-1}b_p}$ and $\Delta_1=D_{n_3}+\Span{\Ad_a^{n_3}b_p}$ we obtain $\Delta_0=\Span{b_2}$ and $\Delta_1=\Span{b_1,b_2,[a,b_2]}$, respectively. With these distributions, the items \ref{it:a} to \ref{it:e} of Theorem \ref{thm:1} are met. We have
			\begin{align*}
				\begin{aligned}
					\Delta_1^{(1)}&=\Span{\epsilon\cos(\theta)\partial_x+\epsilon\sin(\theta)\partial_z+\partial_\theta,\partial_{v_x},\partial_{v_z},\partial_\omega}=\overline{\Delta}_1\,,
				\end{aligned}
			\end{align*}
			thus, item \ref{it:b} is met and we have $n_2=3$. Furthermore, we have $G_1=\overline{\Delta}_1+[a,\overline{\Delta}_1]=\mathcal{T}(\mathcal{X})$, therefore $\Dim{G_1}=\Dim{\overline{\Delta}_1}+2$ holds and thus, the $x_1$-subsystem in a corresponding triangular form \eqref{eq:triangularFormDifferentLengths} consists of two integrator chains, $G_1=\mathcal{T}(\mathcal{X})$ furthermore implies that both of these chains are of length one. Thus, according to Section \ref{se:determinifFlatOutputs}, Case 1, flat outputs compatible with the triangular form are all pairs of functions $(\varphi^1,\varphi^2)$, which satisfy $\Span{\D\varphi^1,\D\varphi^2}=(\overline{\Delta}_1)^\perp$. From $(\overline{\Delta}_1)^\perp=\Span{\epsilon\sin(\theta)\D\theta-\D z,\epsilon\cos(\theta)\D\theta-\D x}$, a possible pair of such functions follows as \eg $\varphi^1=\epsilon\cos(\theta)+z$, $\varphi^2=\epsilon\sin(\theta)-x$. In conclusion, the planar VTOL \eqref{eq:vtol} is static feedback equivalent to the triangular form \eqref{eq:triangularFormDifferentLengths} with $n_{1,1}=n_{1,2}=1$, $n_2=3$ and $n_3=1$. Let us demonstrate the transformation of the VTOL \eqref{eq:vtol} into the form \eqref{eq:triangularFormDifferentLengths}, such that the components of the flat output $\varphi^1=\epsilon\cos(\theta)+z$, $\varphi^2=\epsilon\sin(\theta)-x$ appear as top variables, by following the six steps of the sufficiency part of the proof of Theorem \ref{thm:1}. 
			\subparagraph*{Step 1:} In this example, the sequence of involutive distributions \eqref{eq:involutiveSequenceProof} reduces to
			\begin{align}\label{eq:distributionsVtol}
				\begin{aligned}
					\Delta_0\subset\overline{\Delta}_1\subset G_1=\mathcal{T}(\mathcal{X})\,.
				\end{aligned}
			\end{align}
			These distributions can be straightened out by the state transformation 
			\begin{align}\label{eq:vtolTransformation1}
				\begin{aligned}
					x_{1,1}^1&=\varphi^1=\epsilon\cos(\theta)+z&&&x_2^2&=\Lie_a\varphi^2=\epsilon\omega\cos(\theta)-v_x\\
					x_{1,2}^1&=\varphi^2=\epsilon\sin(\theta)-x&&&x_2^3&=\theta\\
					x_2^1&=\Lie_a\varphi^1=-\epsilon\omega\sin(\theta)+v_z&&&x_3^1&=\omega\,.
				\end{aligned}
			\end{align}
			\begin{remark}
				Choosing $x_{1,1}^1=\varphi^1$, $x_{1,2}^1=\varphi^2$, $x_2^1=\Lie_a\varphi^1$ and $x_2^2=\Lie_a\varphi^2$, is not mandatory for straightening out the distributions \eqref{eq:distributionsVtol}. It is a short cut which immediately transforms the $x_1$-subsystem into Brunovsky normal form, \ie it joins together the Steps 1 to 3 of the proof. Alternatively, we could choose a transformation which just straightens out the distributions \eqref{eq:distributionsVtol}, then introduce the components of the flat output as the states of the $x_1$-subsystem and then apply the Steps 2 and 3.
			\end{remark}\noindent
			Applying the transformation \eqref{eq:vtolTransformation1} to \eqref{eq:vtol} results in
			\begin{align}\label{eq:vtolStep3}
				\begin{aligned}
					&f_1:\quad\begin{aligned}
						\dot{x}_{1,1}^1&=x_2^1\\
						\dot{x}_{1,2}^1&=x_2^2
					\end{aligned}&&\qquad f_2:\quad\begin{aligned}
						\dot{x}_2^1&=\cos(x_2^3)(u^1-\epsilon(x_3^1)^2)-1\\
						\dot{x}_2^2&=\sin(x_2^3)(u^1-\epsilon(x_3^1)^2)\\
						\dot{x}_2^3&=x_3^1
					\end{aligned}&&\qquad f_3:\quad\begin{aligned}
						\dot{x}_3^1&=u^2\,.
					\end{aligned}
				\end{aligned}
			\end{align}
			\subparagraph*{Step 4:} Since we have $n_3=1$, in order to normalize the first equation of the $x_2$-subsystem, we have to apply the input transformation explained in footnote \ref{fn:transformation}, instead of the state transformation \eqref{eq:firstInput}. The input transformation reads $\tilde{u}^2=\cos(x_2^3)(u^1-\epsilon(x_3^1)^2)-1$, $\tilde{u}^1=u^2$. Applying this transformation to \eqref{eq:vtolStep3} results in
			\begin{align*}
				&f_1:\quad\begin{aligned}
					\dot{x}_{1,1}^1&=x_2^1\\
					\dot{x}_{1,2}^1&=x_2^2
				\end{aligned}&&\qquad f_2:\quad\begin{aligned}
					\dot{x}_2^1&=\tilde{u}^2\\
					\dot{x}_2^2&=\tan(x_2^3)(\tilde{u}^2+1)\\
					\dot{x}_2^3&=x_3^1
				\end{aligned}&&\qquad f_3:\quad\begin{aligned}
					\dot{x}_3^1&=\tilde{u}^1\,.
				\end{aligned}
			\end{align*}
			\subparagraph*{Step 5:} We have to successively introduce the components of the input vector field associated with the input $\tilde{u}^2$ of the $x_2$-subsystem as new states. Since $n_2=3$, in this example we have only one such step, namely introducing $\tilde{x}_2^3=\tan(x_2^3)$, which results in
			\begin{align*}
				&f_1:\quad\begin{aligned}
					\dot{x}_{1,1}^1&=x_2^1\\
					\dot{x}_{1,2}^1&=x_2^2
				\end{aligned}&&\qquad f_2:\quad\begin{aligned}
					\dot{x}_2^1&=\tilde{u}^2\\
					\dot{x}_2^2&=\tilde{x}_2^3(\tilde{u}^2+1)\\
					\dot{\tilde{x}}_2^3&=(1+(\tilde{x}_2^3)^2)x_3^1
				\end{aligned}&&\qquad f_3:\quad\begin{aligned}
					\dot{x}_3^1&=\tilde{u}^1\,.
				\end{aligned}
			\end{align*}
			The last equation of the $x_2$-subsystem is normalized by introducing $x_{3,1}^1=(1+(\tilde{x}_2^3)^2)x_3^1$, which results in
			\begin{align*}
				&f_1:\quad\begin{aligned}
					\dot{x}_{1,1}^1&=x_2^1\\
					\dot{x}_{1,2}^1&=x_2^2
				\end{aligned}&&\qquad f_2:\quad\begin{aligned}
					\dot{x}_2^1&=\tilde{u}^2\\
					\dot{x}_2^2&=\tilde{x}_2^3(\tilde{u}^2+1)\\
					\dot{\tilde{x}}_2^3&=x_{3,1}^1
				\end{aligned}&&\qquad f_3:\quad\begin{aligned}
					\dot{x}_3^1&=2\tfrac{\tilde{x}_2^3}{1+(\tilde{x}_2^3)^2}(x_{3,1}^1)^2+(1+(\tilde{x}_2^3)^2)\tilde{u}^1\,.
				\end{aligned}
			\end{align*}
			\subparagraph*{Step 6:} By introducing $\tilde{\tilde{u}}^1=2\tfrac{\tilde{x}_2^3}{1+(\tilde{x}_2^3)^2}(x_{3,1}^1)^2+(1+(\tilde{x}_2^3)^2)\tilde{u}^1$, the $x_3$-subsystem takes Brunovsky normal form, \ie $x_{3,1}^1=\tilde{\tilde{u}}^1$, which completes the transformation into the form \eqref{eq:triangularFormDifferentLengths}. The individual transformation steps summarized to one state and input transformation read
			\begin{align}\label{eq:tges2}
				\begin{aligned}
					x_{1,1}^1&=\epsilon\cos(\theta)+z&&&	\tilde{x}_2^3&=\tan(\theta)\\
					x_{1,2}^1&=\epsilon\sin(\theta)-x&&&x_{3,1}^1&=\tfrac{\omega}{\cos^2(\theta)}\\
					x_2^1&=-\epsilon\omega\sin(\theta)+v_z&&&\tilde{\tilde{u}}^1&=\tfrac{2\omega^2\sin(\theta)+\cos(\theta)u^2}{\cos^3(\theta)}\\
					x_2^2&=\epsilon\omega\cos(\theta)-v_x&&&\tilde{u}^2&=(u^1-\epsilon\omega^2)\cos(\theta)-1\,.
				\end{aligned}
			\end{align}
		\subsection{Academic example}
			Consider the system
			\begin{align}\label{eq:sin}
				\begin{aligned}
					\dot{x}^1&=u^1\\
					\dot{x}^2&=u^2\\
					\dot{x}^3&=\sin(\tfrac{u^1}{u^2})\,,
				\end{aligned}
			\end{align}
			also considered in \cite{Levine:2009} and \cite{Schoberl:2014}. This system is not an AI-system, however, every non-linear system of the general form $\dot{x}=f(x,u)$ becomes an AI-system with according properties regarding flatness, by one-fold prolonging every control. By prolonging both controls of \eqref{eq:sin}, we obtain the AI-system
			\begin{align}\label{eq:sin_ai}
				\begin{aligned}
					\dot{x}^1&=u^1\\
					\dot{x}^2&=u^2\\
					\dot{x}^3&=\sin(\tfrac{u^1}{u^2})\\
					\dot{u}^1&=u^1_1\\
					\dot{u}^2&=u^2_1\,,
				\end{aligned}
			\end{align}
			with the new state $z=(x^1,x^2,x^3,u^1,u^2)$, the new inputs $u^1_1$ and $u^2_1$, the input vector fields $b_1=\partial_{u^1}$ and $b_2=\partial_{u^2}$ and the drift $a=u^1\partial_{x^1}+u^2\partial_{x^2}+\sin(\tfrac{u^1}{u^2})\partial_{x^3}$. In the following, we show that \eqref{eq:sin_ai} is static feedback equivalent to the triangular form \eqref{eq:triangularFormDifferentLengths} by applying Theorem \ref{thm:1}. The distribution $D_1=\Span{b_1,b_2}$ is involutive, the distribution $D_2=\Span{\partial_{u^1},\partial_{u^2},u^2\partial_{x^1}+\cos(\tfrac{u^1}{u^2})\partial_{x^3},u^1\partial_{x^1}+u^2\partial_{x^2}}$ is not involutive, so we have $n_3=1$ and the conditions $\Dim{D_i}=2i$, $i=1,2$ are met. The condition $\C{D_2}\neq D_1$ of Theorem \ref{thm:1} is also met. Before we can evaluate the remaining conditions of Theorem \ref{thm:1}, we have to construct a vector field $b_p$ for this system. For the distribution $H=D_2+[D_1,D_2]$, see Remark \ref{rem:simplerMethod}, we obtain $H=\mathcal{T}(\mathcal{X})$. Thus, a vector field $b_p$ cannot be determined via the method described in Remark \ref{rem:simplerMethod}. Instead, we determine a vector field $b_p$ via \eqref{eq:alg}. In this example, we have
			\begin{align*}
				\begin{aligned}
					v_1&=\Ad_a^{n_3-1}b_1=b_1=\partial_{u^1}\\
					v_2&=\Ad_a^{n_3-1}b_2=b_2=\partial_{u^2}\,.
				\end{aligned}
			\end{align*}
			By inserting those vector fields into \eqref{eq:alg}, we obtain
			\begin{align}\label{eq:eqsys_sin}
				\begin{aligned}
					&\left((\alpha^1)^2\sin(\tfrac{u^1}{u^2})(u^2)^2+2\alpha^1\alpha^2(\cos(\tfrac{u^1}{u^2})u^2-\sin(\tfrac{u^1}{u^2})u^1)u^2+\right.\\
					&\hspace{8em}\left.(\alpha^2)^2(\sin(\tfrac{u^1}{u^2})u^1-2\cos(\tfrac{u^1}{u^2})u^2)u^1\right)\partial_{x^3}\overset{!}{=}0\Mod D_{2}\,.
				\end{aligned}
			\end{align}
			The condition \eqref{eq:eqsys_sin} admits two independent non-trivial solutions, namely $\alpha^1=\lambda u^1$, $\alpha^2=\lambda u^2$ and $\alpha^1=\lambda u^1\tan(\tfrac{u^1}{u^2})-2u^2$, $\alpha^2=\lambda u^2\tan(\tfrac{u^1}{u^2})$, both solutions with an arbitrary non-zero function $\lambda(z)$. Thus, with the choice $\lambda=1$, we obtain the candidates $b_p=u^1\partial_{u^1}+u^2\partial_{u^2}$ and $b_p=(u^1\tan(\tfrac{u^1}{u^2})-2u^2)\partial_{u^1}+u^2\tan(\tfrac{u^1}{u^2})\partial_{u^2}$. If \eqref{eq:sin_ai} is indeed static feedback equivalent to \eqref{eq:triangularFormDifferentLengths}, the remaining conditions of Theorem \ref{thm:1} must be met with at least one of these candidates. This is indeed the case, namely with the vector field $b_p=u^1\partial_{u^1}+u^2\partial_{u^2}$, \ie the vector field constructed from the first solution. For the distributions $\Delta_0=D_{n_3-1}+\Span{\Ad_a^{n_3-1}b_p}$ and $\Delta_1=D_{n_3}+\Span{\Ad_a^{n_3}b_p}$ we obtain
			\begin{align*}
				\begin{aligned}
					\Delta_0&=\Span{b_p}=\Span{u^1\partial_{u^1}+u^2\partial_{u^2}}\\
					\Delta_1&=\Span{b_1,b_2,[a,b_p]}=\Span{\partial_{u^1},\partial_{u^2},u^1\partial_{x^1}+u^2\partial_{x^2}}\,.
				\end{aligned}
			\end{align*}
			With these distributions, the items \ref{it:a} to \ref{it:e} of Theorem \ref{thm:1} are met. We have
			\begin{align*}
				\begin{aligned}
					\Delta_1^{(1)}&=\Span{\partial_{u^1},\partial_{u^2},\partial_{x^1},\partial_{x^2}}=\overline{\Delta}_1\,,
				\end{aligned}
			\end{align*}
			thus, item \ref{it:b} is met and we have $n_2=3$. Furthermore, we have $G_1=\overline{\Delta}_1+[a,\overline{\Delta}_1]=\mathcal{T}(\mathcal{X})$, therefore $\Dim{G_1}=\Dim{\overline{\Delta}_1}+1$ holds and thus, the $x_1$-subsystem in a corresponding triangular form \eqref{eq:triangularFormDifferentLengths} only consists of one integrator chain, $G_1=\mathcal{T}(\mathcal{X})$, \ie $s=1$, furthermore implies that this chain is of length one. Thus, according to Section \ref{se:determinifFlatOutputs}, Case 2, flat outputs compatible with the triangular form  are all pairs of functions $(\varphi^1,\varphi^2)$, which satisfy $L^\perp=\Span{\D\varphi^1,\D\Lie_a\varphi^1,\D\varphi^2}=\Delta_1^\perp+\Span{\D\Lie_a\varphi^1}$ with $\Span{\D\varphi^1}=(\overline{\Delta}_1)^\perp$. We have $(\overline{\Delta}_1)^\perp=\Span{\D x^3}$, thus $\varphi^1=\varphi^1(x^3)$. Furthermore, we have $\Lie_a\varphi^1(x^3)=\sin(\tfrac{u^1}{u^2})\partial_{x^3}\varphi^1(x^3)$ and thus $L^\perp=\Delta_1^\perp+\Span{\D\Lie_a\varphi^1}=\Span{u^2\D x^1-u^1\D x^2,u^2\D u^1-u^1\D u^2,\D x^3}$. Therefore, $\varphi^2=\varphi^2(u^1/u^2,x^1-x^2u^1/u^2,x^3)$, chosen such that $\D\varphi^1\wedge\D\Lie_a\varphi^1\wedge\D\varphi^2\neq 0$. A possible flat output is thus \eg $\varphi^1=x^3$, $\varphi^2=x^1-x^2u^1/u^2$. In conclusion, \eqref{eq:sin_ai} is static feedback equivalent to the triangular form \eqref{eq:triangularFormDifferentLengths} with $n_1=1$, $n_2=3$ and $n_3=1$. Indeed, by applying a suitable state and input transformation to \eqref{eq:sin_ai}, which again can be derived systematically following the six steps of the sufficiency part of the proof of Theorem \ref{thm:1}, the system \eqref{eq:sin_ai} takes the form
			\begin{align*}
				\begin{aligned}
					&f_1:\quad\begin{aligned}
						\dot{z}_{1,1}^1&=z_3^1\\
					\end{aligned}&&\qquad f_2:\quad\begin{aligned}
						\dot{z}_2^1&=\tilde{u}^2_1\\
						\dot{z}_2^2&=z_2^3\tilde{u}^2_1\\
						\dot{z}_2^3&=z_{3,1}^1+\tfrac{z_2^1z_2^3}{1-(z_2^1)^2}\tilde{u}^2_1
					\end{aligned}&&\qquad f_3:\quad\begin{aligned}
						\dot{z}_{3,1}^1&=\tilde{u}^1_1\,,
					\end{aligned}
				\end{aligned}
			\end{align*}
			which is of the form \eqref{eq:triangularFormDifferentLengths}.
		\subsection{Further academic examples}
			Consider the following two academic examples
			\begin{align*}
				\begin{aligned}
					\dot{x}^1&=u^1\\
					\dot{x}^2&=u^2\\
					\dot{x}^3&=u^1u^2
				\end{aligned}&&&
				\begin{aligned}
					\dot{x}^1&=u^1\\
					\dot{x}^2&=u^2\\
					\dot{x}^3&=\sqrt{u^1u^2}\,,
				\end{aligned}
			\end{align*}
			which are similar to \eqref{eq:sin} in the previous section and are also treated in \eg \cite{Schoberl:2014}. Also these systems are static feedback equivalent to \eqref{eq:triangularFormDifferentLengths} (after turning them into AI-systems by prolonging each of their controls, as demonstrated on the previous example) and thus, can be transformed into the form \eqref{eq:triangularFormDifferentLengths} systematically. For these systems, the dimensions of the individual subsystems in a corresponding triangular form \eqref{eq:triangularFormDifferentLengths} would be $n_1=1$, $n_2=3$, $n_3=1$ and $n_1=0$, $n_2=4$, $n_3=1$, respectively. Therefore, these systems become static feedback linearizable by prolonging a suitably chose control two-fold (as $n_2=3$) or three-fold (as $n_2=4$), respectively. Flat outputs for these systems can again be derived systematically as described in Section \ref{se:determinifFlatOutputs}, without actually transforming the systems into the form \eqref{eq:triangularFormDifferentLengths}.
		\subsection{Explicit transformation into the triangular form}
			Based on the following academic example, we once more demonstrate the transformation into the triangular form \eqref{eq:triangularFormDifferentLengths} by following the six steps of the sufficiency part of the proof of Theorem \ref{thm:1}. Consider the system
			\begin{align}\label{eq:academic}
				\begin{aligned}
					\dot{x}^1&=x^2&&&\dot{x}^6&=x^7(x^9-x^8x^{10})\\
					\dot{x}^2&=x^4+\sin(x^6)&&&\dot{x}^7&=x^1(x^8x^{10}-x^9)+\sin(x^8)\\
					\dot{x}^3&=x^2+x^5&&&\dot{x}^8&=x^9+x^{10}\\
					\dot{x}^4&=(x^9-x^8x^{10})(1-\cos(x^6)x^7)&&&\dot{x}^9&=u^1\\
					\dot{x}^5&=x^6(x^9-x^8x^{10})&&&\dot{x}^{10}&=u^2\,.
				\end{aligned}
			\end{align}
			The input vector fields are given by $b_1=\partial_{x^9}$ and $b_2=\partial_{x^{10}}$, the drift is given by
			\begin{align*}
				\begin{aligned}
					a&=x^2\partial_{x^1}+(x^4+\sin(x^6))\partial_{x^2}+(x^2+x^5)\partial_{x^3}+(x^9-x^8x^{10})(1-\cos(x^6)x^7)\partial_{x^4}+\\
					&\hspace{2em}x^6(x^9-x^8x^{10})\partial_{x^5}+x^7(x^9-x^8x^{10})\partial_{x^6}+(x^1(x^8x^{10}-x^9)+\sin(x^8))\partial_{x^7}+(x^9+x^{10})\partial_{x^8}\,.
				\end{aligned}
			\end{align*}
			The distributions
			\begin{align*}
				\begin{aligned}
					D_1&=\Span{b_1,b_2}=\Span{\partial_{x^{10}},\partial_{x^{9}}}\\
					D_2&=D_1+[a,D_1]=\Span{\partial_{x^{10}},\partial_{x^{9}},\partial_{x^8},(x^7\cos(x^6)-1)\partial_{x^4}-x^6\partial_{x^5}-x^7\partial_{x^6}+x^1\partial_{x^7}}
				\end{aligned}
			\end{align*}
			are involutive, the distribution
			\begin{align*}
				\begin{aligned}
					D_3&=D_2+[a,D_2]=\Span{\partial_{x^{10}},\partial_{x^{9}},\partial_{x^8},\partial_{x^7},\partial_{x^2}+x^6\partial_{x^3}+\cos(x^6)\sin(x^8)\partial_{x^4}-\sin(x^8)\partial_{x^6},\\
						&\hspace{12em}x^7\partial_{x^2}+x^6x^7\partial_{x^3}+\sin(x^8)\partial_{x^4}+x^6\sin(x^8)\partial_{x^5}}\,,
				\end{aligned}
			\end{align*}
			is not invoultive, so we have $n_3=2$. The conditions $\Dim{D_i}=2i$, $i=1,\ldots,3$ and $\C{D_3}\neq D_2$ are met. For the distribution $H=D_3+[D_2,D_3]$, see Remark \ref{rem:simplerMethod}, we obtain
			\begin{align*}
				\begin{aligned}
					H&=\Span{\partial_{x^{10}},\partial_{x^{9}},\partial_{x^8},\partial_{x^7},x^7\cos(x^6)\partial_{x^2}+\sin(x^8)\partial_{x^6},x^7\partial_{x^2}-x^6\sin(x^8)\partial_{x^5},\partial_{x^2}+x^6\partial_{x^3},\\
						&\hspace{6em}x^7\partial_{x^2}+\sin(x^8)\partial_{x^4}}\,.
				\end{aligned}
			\end{align*}
			For the vector fields $\Ad_a^{n_3+1}b_1$ and $\Ad_a^{n_3+1}b_2$ we have
			\begin{align*}
				\begin{aligned}
					\Ad_a^3b_1&=-\partial_{x^1}+\tfrac{1}{x^6}\partial_{x^2}\Mod H\\
					\Ad_a^3b_2&=x^8\partial_{x^1}-\tfrac{x^8}{x^6}\partial_{x^2}\Mod H\,.
				\end{aligned}
			\end{align*}
			The linear combination $x^8\Ad_a^2b_1+\Ad_a^2b_2=0\Mod H$ is obviously contained in $H$. Thus, we have $b_p=x^8b_1+b_2=x^8\partial_{x^9}+\partial_{x^{10}}$ and for the distributions $\Delta_0=D_{n_3-1}+\Span{\Ad_a^{n_3-1}b_p}$ and $\Delta_1=D_{n_3}+\Span{\Ad_a^{n_3}b_p}$ we obtain
			\begin{align*}
				\begin{aligned}
					\Delta_0&=\Span{\partial_{x^{10}},\partial_{x^{9}},\partial_{x^8}}\\
					\Delta_1&=\Span{\partial_{x^{10}},\partial_{x^{9}},\partial_{x^8},\partial_{x^7},(1-x^7\cos(x^6))\partial_{x^4}+x^6\partial_{x^5}+x^7\partial_{x^6}}\,.
				\end{aligned}
			\end{align*}
			With these distributions the items \ref{it:a} to \ref{it:e} of Theorem \ref{thm:1} are met. We have
			\begin{align*}
				\begin{aligned}
					\Delta_1^{(1)}&=\Span{\partial_{x^{10}},\partial_{x^{9}},\partial_{x^8},\partial_{x^7},\cos(x^6)\partial_{x^4}-\partial_{x^6},\partial_{x^4}+x^6\partial_{x^5}}\\
					\Delta_1^{(2)}&=\Span{\partial_{x^{10}},\ldots,\partial_{x^4}}=\overline{\Delta}_1\,,
				\end{aligned}
			\end{align*}
			thus, item \ref{it:b} is met and we have $n_2=4$. Furthermore, we have $G_1=\Span{\partial_{x^{10}},\ldots,\partial_{x^2}}$, $G_2=\mathcal{T}(\mathcal{X})$ and $\Dim{G_1}=\Dim{\overline{\Delta}_1}+2$ holds. Thus, in a corresponding triangular form \eqref{eq:triangularFormDifferentLengths}, the $x_1$-subsystem consists of two integrator chains with the lengths one and two. Thus, according to Section \ref{se:determinifFlatOutputs}, Case 1, flat outputs compatible with the triangular form are all pairs of functions $(\varphi^1,\varphi^2)$, which satisfy $\Span{\D\varphi^1}=G_1^\perp$ and $\Span{\D\varphi^1,\D\Lie_a\varphi^1,\D\varphi^2}=(\overline{\Delta}_1)^\perp$. From $G_1^\perp=\Span{\D x^1}$, $\varphi^1=\varphi^1(x^1)$ follows. From $\Lie_a\varphi^1(x^1)=x^2\partial_{x^1}\varphi^1(x^1)$ and  $(\overline{\Delta}_1)^\perp=\Span{\D x^1,\D x^2,\D x^3}$, it follows that $\varphi^2=\varphi^2(x^1,x^2,x^3)$, chosen such that $\D\varphi^1\wedge\D\Lie_a\varphi^1\wedge\varphi^2\neq 0$. A possible flat output is thus \eg $\varphi^1=x^1$, $\varphi^2=x^3$. In the following, we transform \eqref{eq:academic} into the triangular form \eqref{eq:triangularFormDifferentLengths}, such that the components of the flat output $\varphi^1=x^1$, $\varphi^2=x^3$ appear as top variables in the triangular form.
			\subparagraph*{Step 1:} In this example, the distributions
			\begin{align*}
				\begin{aligned}
					D_1\subset\Delta_0\subset\C{\Delta_1^{(1)}}\subset\overline{\Delta}_1\subset G_1\subset G_2=\mathcal{T}(\mathcal{X})\,,
				\end{aligned}
			\end{align*}
			corresponding to the sequence \eqref{eq:involutiveSequenceProof}, are already straightened out. Therefore, \eqref{eq:academic} is structurally already in the form \eqref{eq:step1}. Indeed, by renaming the states according to
			\begin{align*}
				\begin{aligned}
					\begin{aligned}
						x_1^1&=x_1\\
						x_1^2&=x_2\\
						x_1^3&=x_3
					\end{aligned}
					&&\quad\begin{aligned}
						x_2^1&=x^4\\
						x_2^2&=x^5\\
						x_2^3&=x^6\\
						x_2^4&=x^7
					\end{aligned}
					&&\quad\begin{aligned}
						x_3^1&=x^8\\
						x_3^2&=x^9\\
						x_3^3&=x^{10}\,,
					\end{aligned}
				\end{aligned}
			\end{align*}
			we obtain
			\begin{align*}
				\begin{aligned}
					f_1:&\quad\begin{aligned}
						\dot{x}_1^1&=x_1^2\\
						\dot{x}_1^2&=x_2^1+\sin(x_2^3)\\
						\dot{x}_1^3&=x_1^2+x_2^2
					\end{aligned}\\[1ex]
					f_2:&\quad\begin{aligned}
						\dot{x}_2^1&=(x_3^2-x_3^1x_3^3)(1-\cos(x_2^3)x_2^4)\\
						\dot{x}_2^2&=x_2^3(x_3^2-x_3^1x_3^3)\\
						\dot{x}_2^3&=x_2^4(x_3^2-x_3^1x_3^3)\\
						\dot{x}_2^4&=x_1^1(x_3^1x_3^3-x_3^2)+\sin(x_3^1)
					\end{aligned}\\[1ex]
					f_3:&\quad\begin{aligned}
						\dot{x}_3^1&=x_3^2+x_3^3\\
						\dot{x}_3^2&=u^1\\
						\dot{x}_3^3&=u^2\,,
					\end{aligned}
				\end{aligned}
			\end{align*}
			which is exactly the form \eqref{eq:step1}. The rank conditions  $\Rank{\partial_{(x_2^1,x_2^2,x_2^3)}f_1}=2$, $\Rank{\partial_{(x_3^1,x_3^2,x_3^3)}f_2}=2$ and $\Rank{\partial_{(x_3^2,x_3^3)}f_2}=1$ hold.
			\subparagraph*{Step 2:} The $x_1$-subsystem is already in Brunovsky normal form except for a normalization of the "inputs" of the integrator chains. To obtain exactly the representation \eqref{eq:step2}, we only have to rename the states of the $x_1$-subsystem according to $x_{1,1}^1=x_1^1$, $x_{1,1}^2=x_1^2$ and $x_{1,2}^1=x_1^3$. This results in
			\begin{align*}
				\begin{aligned}
					f_1:\quad\begin{aligned}
						\dot{x}_{1,1}^1&=x_{1,1}^2\\
						\dot{x}_{1,1}^2&=x_2^1+\sin(x_2^3)\\
						\dot{x}_{1,2}^1&=x_{1,1}^2+x_2^2\,.
					\end{aligned}
				\end{aligned}
			\end{align*}		
			\subparagraph*{Step 3:} The transformation of the $x_1$-subsystem into Brunovsky normal form is completed by normalizing the last two equations of the $x_1$-subsystem, \ie by introducing $\tilde{x}_2^1=x_2^1+\sin(x_2^3)$ and $\tilde{x}_2^2=x_{1,1}^2+x_2^2$, resulting in
			\begin{align*}
				\begin{aligned}
					f_1:&\quad\begin{aligned}
						\dot{x}_{1,1}^1&=x_{1,1}^2\\
						\dot{x}_{1,1}^2&=\tilde{x}_2^1&&&\dot{x}_{1,2}^1&=\tilde{x}_2^2
					\end{aligned}\\[1ex]
					f_2:&\quad\begin{aligned}
						\dot{\tilde{x}}_2^1&=x_3^2-x_3^1x_3^3\\
						\dot{\tilde{x}}_2^2&=x_2^3(x_3^2-x_3^1x_3^3)+\tilde{x}_2^1\\
						\dot{x}_2^3&=x_2^4(x_3^2-x_3^2x_3^3)\\
						\dot{x}_2^4&=x_1^1(x_3^1x_3^3-x_3^2)+\sin(x_3^1)\,.
					\end{aligned}
				\end{aligned}
			\end{align*}
			\subparagraph*{Step 4:} Next, we normalize the first equation of the $x_2$-subsystem, by introducing $x_{3,2}^1=x_3^2-x_3^1x_3^3$. This leads to
			\begin{align*}
				\begin{aligned}
					f_2:&\quad\begin{aligned}
						\dot{\tilde{x}}_2^1&=x_{3,2}^1\\
						\dot{\tilde{x}}_2^2&=x_2^3x_{3,2}^1+\tilde{x}_2^1\\
						\dot{x}_2^3&=x_2^4x_{3,2}^1\\
						\dot{x}_2^4&=\sin(x_3^1)-x_{1,1}^1x_{3,2}^1\,.
					\end{aligned}\\[1ex]
					f_3:&\quad\begin{aligned}
						\dot{x}_3^1&=(1+x_3^1)x_3^3+x_{3,2}^1\\
						\dot{x}_{3,2}^1&=h(x_3^1,x_{3,2}^1,x_3^3,u^1,u^2)\\
						\dot{x}_3^3&=u^2\,,
					\end{aligned}
				\end{aligned}
			\end{align*}
			the $x_2$-subsystem is indeed of the form \eqref{eq:subsys2triangular}.
			\subparagraph*{Step 5:} We have to successively introduce the components of the input vector field associated with the input $x_{3,2}^1$ of the $x_2$-subsystem as new states (which here is actually already the case, we only have to rename the states according to $\tilde{x}_2^3=x_2^3$ and $\tilde{x}_2^4=x_2^4$ to be consistent with the notation in the proof of Theorem \ref{thm:1}). Normalizing the last equation of the $x_2$-subsystem, \ie introducing $x_{3,1}^1=\sin(x_3^1)-x_{1,1}^1x_{3,2}^1$, would complete the transformation of the $x_2$-subsystem to extended chained form. However, this transformation would result in
			\begin{align*}
				\begin{aligned}
					f_2:&\quad\begin{aligned}
						\dot{\tilde{x}}_2^1&=x_{3,2}^1\\
						\dot{\tilde{x}}_2^2&=\tilde{x}_2^3x_{3,2}^1+\tilde{x}_2^1\\
						\dot{\tilde{x}}_2^3&=\tilde{x}_2^4x_{3,2}^1\\
						\dot{\tilde{x}}_2^4&=x_{3,1}^1\,.
					\end{aligned}\\[1ex]
					f_3:&\quad\begin{aligned}
						\dot{x}_{3,1}^1&=h^1(x_{1,1}^1,x_{1,1}^2,x_{3,1}^1,x_{3,2}^1,x_3^3,u^1,u^2)\\
						\dot{x}_{3,2}^1&=h^2(x_{1,1}^1,x_{3,1}^1,x_{3,2}^1,x_3^3,u^1,u^2)\\
						\dot{x}_3^3&=u^2\,,
					\end{aligned}
				\end{aligned}
			\end{align*}
			preventing us from transforming the $x_3$-subsystem into Brunovsky normal from by successively introducing new coordinates from top to bottom, since the inputs $u^1$ and $u^2$ occur in all three equations of the $x_3$-subsystem (the distribution $D_1$ is not straightened out anymore, see also Remark \ref{rem:complete}). Instead, we only introduce $x_{3,1}^1=\sin(x_3^1)$, which results in
			\begin{align*}
				\begin{aligned}
					f_2:&\quad\begin{aligned}
						\dot{\tilde{x}}_2^1&=x_{3,2}^1\\
						\dot{\tilde{x}}_2^2&=\tilde{x}_2^3x_{3,2}^1+\tilde{x}_2^1\\
						\dot{\tilde{x}}_2^3&=\tilde{x}_2^4x_{3,2}^1\\
						\dot{\tilde{x}}_2^4&=x_{3,1}^1-x_{1,1}^1x_{3,2}^1\,.
					\end{aligned}\\[1ex]
					f_3:&\quad\begin{aligned}
						\dot{x}_{3,1}^1&=\sqrt{1-(x_{3,1}^1)^2}\left(1+\arcsin(x_{3,1}^1)\right)x_3^3+x_{3,2}^1\\
						\dot{x}_{3,2}^1&=-x_3^3\left(x_{3,2}^1+x_3^3\left(1+\arcsin(x_{3,1}^1)\right)\right)-\arcsin(x_{3,1}^1)u^{2}+u^{1}\\
						\dot{x}_3^3&=u^2
					\end{aligned}
				\end{aligned}
			\end{align*}
			and keeps $D_1$ straightened out, so the inputs $u^1$ and $u^2$ still only occur in the last two equations of the $x_3$-subsystem.
			\subparagraph*{Step 6:} The last step is to transform the $x_3$-subsystem into Brunovsky normal form. For that, we first introduce $x_{3,1}^2=\sqrt{1-(x_{3,1}^1)^2}\left(1+\arcsin(x_{3,1}^1)\right)x_3^3+x_{3,2}^1$, to obtain
			\begin{align*}
				\begin{aligned}
					f_3:&\quad\begin{aligned}
						\dot{x}_{3,1}^1&=x_{3,1}^2\\
						\dot{x}_{3,1}^2&=h^1(x_{3,1}^1,x_{3,1}^2,u^1,u^2)\\
						\dot{x}_{3,2}^1&=h^2(x_{3,1}^1,x_{3,1}^2,x_{3,2}^1,u^1,u^2)\,.
					\end{aligned}
				\end{aligned}
			\end{align*}
			Finally, we complete the transformation by introducing $\tilde{u}^1=h^1$ and $\tilde{u}^2=h^2$. After applying this input transformation, the complete system reads
			\begin{align*}
				\begin{aligned}
					f_1:&\quad\begin{aligned}
						\dot{x}_{1,1}^1&=x_{1,1}^2\\
						\dot{x}_{1,1}^2&=\tilde{x}_2^1&&&\dot{x}_{1,2}^1&=\tilde{x}_2^2
						\end{aligned}\\[1ex]
					f_2:&\quad\begin{aligned}
						\dot{\tilde{x}}_2^1&=x_{3,2}^1\\
						\dot{\tilde{x}}_2^2&=\tilde{x}_2^3x_{3,2}^1+\tilde{x}_2^1\\
						\dot{\tilde{x}}_2^3&=\tilde{x}_2^4x_{3,2}^1\\
						\dot{\tilde{x}}_2^4&=x_{3,1}^1-x_{1,1}^1x_{3,2}^1\,.
					\end{aligned}\\[1ex]
					f_3:&\quad\begin{aligned}
						\dot{x}_{3,1}^1&=x_{3,1}^2&&&\dot{x}_{3,2}^1&=\tilde{u}^2\\
						\dot{x}_{3,1}^2&=\tilde{u}^1\,.
					\end{aligned}
				\end{aligned}
			\end{align*}
			which is of the form \eqref{eq:triangularFormDifferentLengths}.
	\bibliographystyle{apacite}
	\bibliography{Bibliography}
	\appendix
	\section{Supplements}\label{ap:additionalProofs}
		In this section, details omitted in the proof of our main theorem and proofs concerning the construction of a vector field $b_p$ are provided.
		\subsection*{Proof of Lemma \ref{lem:characteristic}}
			Let $\{c_1,\ldots,c_{n_c}\}$ be a basis for $\C{D}$ and $\{v_1,\ldots,v_{d}\}$ a basis for $D$. We obviously have $[c_i,v_j]\in D$ and $[v_j,v_k]\in D^{(1)}$. From the Jacobi identity
			\begin{align*}
				\underbrace{[v_j,\underbrace{[c_i,v_k]}_{\in D}]}_{\in D^{(1)}}+\underbrace{[v_k,\underbrace{[v_j,c_i]}_{\in D}]}_{\in D^{(1)}}+[c_i,\underbrace{[v_k,v_j]}_{\in D^{(1)}}]&=0\,,&i\in\{1,\ldots,n_c\}\,,~j,k\in\{1,\ldots,d\}\,,
			\end{align*}
			it follows that $[c_i,[v_k,v_j]]\in D^{(1)}$. Thus, every vector field $c_i\in\C{D}$ is also characteristic for $D^{(1)}$, \ie $\C{D}\subset\C{D^{(1)}}$. 
		\subsection*{Proof of Lemma \ref{lem:cartan}}
			Let us construct a special basis for the distribution $D$, namely $D=\Span{c_1,\ldots,c_{d-2},v_1,v_2}$, with $c_j\in\C{D}$. Because of $\Dim{D^{(i)}}=d+i$, bases for $D^{(1)}$ and $D^{(2)}$ are then given by $D^{(1)}=\Span{c_1,\ldots,c_{d-2},v_1,v_2,v_3}$ and $D^{(2)}=\Span{c_1,\ldots,c_{d-2},v_1,v_2,v_3,v_4}$ (with $v_3=[v_1,v_2]$ and $v_4=[v_1,v_3]$ or $v_4=[v_2,v_3]$ if $[v_1,v_3]=0\Mod D^{(1)}$). We obviously have $[v_1,v_2]=0\Mod D^{(1)}$. Furthermore, we have $[v_1,v_3]=\alpha^1v_4\Mod D^{(1)}$ and $[v_2,v_3]=\alpha^2v_4\Mod D^{(1)}$, where $\alpha^1$ and $\alpha^2$ are some functions and at least $\alpha^1\neq 0$ or $\alpha^2\neq 0$. The vector field $\tilde{v}_1=\alpha^2v_1-\alpha^1v_2\in D$ satisfies $[\tilde{v}_1,v_i]=0\Mod D^{(1)}$, $i=1,\ldots,3$ and $[\tilde{v}_1,c_j]=0\Mod D\subset D^{(1)}$, $j=1,\ldots,d-2$. Thus, $\tilde{v}_1$ is a characteristic vector field of $D^{(1)}$. Because of Lemma \ref{lem:characteristic}, we furthermore have $\C{D}\subset\C{D^{(1)}}$. Thus, we have $\C{D^{(1)}}=\C{D}+\Span{\tilde{v}_1}\subset D$. 
			The rest of the proof follows the same line. A basis for $D^{(1)}$ is given by $D^{(1)}=\Span{c_1,\ldots,c_{d-1},w_1,w_2}$, where $c_{d-1}=\tilde{v}_1$, $w_1=v_1$ (or $w_1=v_2$ if $v_1$ and $\tilde{v}_1$ are collinear$\Mod\C{D}$) and $w_2=v_3$. This way, we formally obtained the same problem as before. Thus, by essentially the same argumentation as before, it follows that $\C{D^{(2)}}=\C{D^{(1)}}+\Span{\tilde{w}_1}\subset D^{(1)}$, with $\tilde{w}_1$ being a suitable linear combination of $w_1$ and $w_2$. Continuing this argumentation, $\Dim{\C{D^{(i)}}}=d-2+i$ and $\C{D^{(i)}}\subset D^{(i-1)}$, $i=1,\ldots,l-1$ follows.
		\subsection{Details omitted in the proof of Theorem \ref{thm:1}}\label{ap:detailsTheorem1}
			\subsection*{Feedback invariance of certain distributions}
			Right before Theorem \ref{thm:1}, for an AI-system \eqref{eq:aisystem}, the distributions $D_i$, $i=1,\ldots,n_3$, $D_1=\Span{b_1,b_2}$ and $D_{i+1}=D_i+[a,D_i]$, with the smallest integer $n_3$ such that $D_{n_3+1}$ is not involutive, are defined. Since they are of importance for Theorem \ref{thm:1}, let us discuss their invariance with respect to invertible static feedback transformations\footnote{The following proof of the feedback invariance of $D_{n_3+1}$ is a replication of a part of the proof of Proposition 7.1 in \cite{NicolauRespondek:2016}, adapted to our notation.}. It is well known that the involutive distributions $D_i$, $i=1,\ldots,n_3$ are feedback invariant, see \eg \cite{JakubczykRespondek:1980} or \cite{NijmeijervanderSchaft:1990}. To show that also the first non-involutive distribution $D_{n_3+1}=D_{n_3}+[a,D_{n_3}]$ is feedback invariant, we calculate $D_{n_3+1}$ after applying an invertible static feedback. After applying an invertible static feedback to an AI-system \eqref{eq:aisystem}, its input vector fields and drift read $\tilde{b}_1=\beta_1^1b_1+\beta_1^2b_2$, $\tilde{b}_2=\beta_2^1b_1+\beta_2^2b_2$, with $\beta_1^1\beta_2^2-\beta_1^2\beta_2^1\neq 0$ and  $\tilde{a}=a+\gamma^1b_1+\gamma^2b_2$, where $\beta_i^j$ and $\gamma^j$ are functions of the state $x$ of the system. Since $\tilde{D}_{n_3}=D_{n_3}$ is involutive, we obtain
			\begin{align*}
				\begin{aligned}
					\tilde{D}_{n_3+1}&=D_{n_3}+[\tilde{a},D_{n_3}]\\
					&=D_{n_3}+[a+\gamma^1b_1+\gamma^2b_2,D_{n_3}]\\
					&=\underbrace{D_{n_3}+[a,D_{n_3}]}_{=D_{n_3+1}}+\underbrace{[\gamma^1b_1+\gamma^2b_2,D_{n_3}]}_{\subset D_{n_3}}\,,
				\end{aligned}
			\end{align*}
			\ie $\tilde{D}_{n_3+1}=D_{n_3+1}$ indeed holds. Also the distributions $\Delta_0=D_{n_3-1}+\Span{\Ad_a^{n_3-1}b_p}$ and $\Delta_1=D_{n_3}+\Span{\Ad_a^{n_3}b_p}$ play a crucial role in Theorem \ref{thm:1}. These are also feedback invariant, \ie for a fixed vector field $b_p\in D_1$, calculating these distribution with the feedback modified drift $\tilde{a}=a+\gamma^1b_1+\gamma^2b_2$, yields the same distributions $\Delta_0$ and $\Delta_1$. Furthermore, only the direction of $b_p$ matters, \ie with $\tilde{b}_p=\lambda b_p$ with an arbitrary non-zero function $\lambda$ of the state of the system, we again obtain the same distributions $\Delta_0$ and $\Delta_1$. To show this, note that because of the involutivity of $D_1,\ldots,D_{n_3}$ and $\Ad_a^{i-1}b_j\in D_i$, we have
			\begin{align*}
				\begin{aligned}
					\Ad_{\tilde{a}}^1\tilde{b}_p&=[a+\gamma^1b_1+\gamma^2b_2,\lambda b_p]=\lambda\Ad_a^1b_p\Mod D_1\\
					\Ad_{\tilde{a}}^2\tilde{b}_p&=[a+\gamma^1b_1+\gamma^2b_2,\lambda\Ad_a^1b_p\Mod D_1]=\lambda\Ad_a^2b_p\Mod D_2\\
					&\vdotswithin{=}\\
					\Ad_{\tilde{a}}^{n_3-1}\tilde{b}_p&=\lambda\Ad_a^{n_3-1}b_p\Mod D_{n_3-1}\\
					\Ad_{\tilde{a}}^{n_3}\tilde{b}_p&=\lambda\Ad_a^{n_3}b_p\Mod D_{n_3}
				\end{aligned}
			\end{align*}
			and thus, the distribution $\tilde{\Delta}_0=D_{n_3-1}+\Span{\Ad_{\tilde{a}}^{n_3-1}\tilde{b}_p}$ indeed coincides with $\Delta_0=D_{n_3-1}+\Span{\Ad_a^{n_3-1}b_p}$, and $\tilde{\Delta}_1=D_{n_3}+\Span{\Ad_{\tilde{a}}^{n_3}\tilde{b}_p}$ indeed coincides with $\Delta_1=D_{n_3}+\Span{\Ad_a^{n_3}b_p}$.
			\subsection*{Details necessity}
				\begin{factProof}[\ref{ft:1}]
					In the necessity part of the proof of Theorem \ref{thm:1}, we claimed that the distribution $D_{n_3+1}=\Span{\partial_{x_3},b_1^c,b_2^c,[a,b_2^c+gb_1^c]}$, with $b_1^c=\partial_{x_2^{n_2}}$ and  $b_2^c=\partial_{x_2^1}+x_2^3\partial_{x_2^2}+\ldots+x_2^{n_2}\partial_{x_2^{n_2-1}}$, is not involutive, meets $\Dim{D_{n_3+1}}=2n_3+2$ and $\C{D_{n_3+1}}\neq D_{n_3}$. To show this, recall that the drift vector field of the complete system \eqref{eq:triangularFormDifferentLengths} is given by $a=f_1+x_{3,2}^1(b_2^c+gb_1^c)+x_{3,1}^1b_1^c+a_2+a_3$ and thus, the vector field $[a,b_2^c+gb_1^c]$ is of the form
					\begin{align*}
						\begin{aligned}
							[a,b_2^c+gb_1^c]&=[f_1,b_2^c]+x_{3,1}^1[b_1^c,b_2^c]+[a_2,b_2^c]+\underbrace{g[a_2,b_1^c]}_{h(x_1,x_2)[b_1^c,b_2^c]}\Mod\Span{\partial_{x_3},b_1^c,b_2^c}\\
							&=[f_1,b_2^c]+\tilde{h}(x_1,x_2,x_{3,1}^1)[b_1^c,b_2^c]+[a_2,b_2^c]\Mod\Span{\partial_{x_3},b_1^c,b_2^c}\,.
						\end{aligned}
					\end{align*}
					Since $[b_1^c,b_2^c]=\partial_{x_2^{n_2-1}}\notin\Span{\partial_{x_3},b_1^c,b_2^c}$, the condition $\Dim{D_{n_3+1}}=2n_3+2$ indeed holds, it holds independently of the dimension $n_1=n_{1,1}+n_{1,2}$ of the $x_1$-subsystem and the actual form of the drift vector field $a_2$ of the $x_2$-subsystem. Next, let us show the non-involutivity of $D_{n_3+1}$ by contradiction. Assume that $D_{n_3+1}$ would be involutive. The distribution $D_{n_3+1}$ contains the vector fields $b_1^c$ and $b_2^c$. Since $[b_1^c,b_2^c]=\partial_{x_2^{n_2-1}}\notin\Span{\partial_{x_3},b_1^c,b_2^c}$, in order for $D_{n_3+1}$ to be involutive, there must hold $[a,b_2^c+gb_1^c]=h(x_1,x_2,x_{3,1}^1)[b_1^c,b_2^c]\Mod\Span{\partial_{x_3},b_1^c,b_2^c}$ (otherwise, the Lie bracket $[b_1^c,b_2^c]$ would not be contained in $D_{n_3+1}$). This is only possible if $[f_1,b_2^c]=0\Mod\Span{\partial_{x_3},b_1^c,b_2^c,[b_1^c,b_2^c]}$ and thus $n_1=0$ (for $n_1>0$, $[f_1,b_2^c]$ has non-zero components in the $\partial_{x_1}$-direction and those are certainly not contained in $\Span{\partial_{x_3},b_1^c,b_2^c,[b_1^c,b_2^c]}$). This shows the non-involutivity of $D_{n_3+1}$ for the case $n_1>0$. For $n_1=0$, the case $[a,b_2^c+gb_1^c]=h(x_1,x_2,x_{3,1}^1)[b_1^c,b_2^c]\Mod\Span{\partial_{x_3},b_1^c,b_2^c}$ can indeed occur (it occurs when $n_1=0$ and $[a_2,b_2^c]=0\Mod\Span{\partial_{x_3},b_1^c,b_2^c,[b_1^c,b_2^c]}$) and results in $D_{n_3+1}=\Span{\partial_{x_3},b_1^c,b_2^c,[b_1^c,b_2^c]}$. This distribution is only involutive, if $[[b_1^c,b_2^c],b_2^c]=\partial_{x_2^{n_2-2}}\in D_{n_3+1}$, which only holds for $n_2=3$ (we do not consider the degenerated case $n_2\leq 2$, see also Remark \ref{rem:dimensions}). However, for $n_2=3$, we have $D_{n_3+1}=\mathcal{T}(\mathcal{X})$ and in turn, the system is static feedback linearizable. In the case that $D_{n_3+1}=\Span{\partial_{x_3},b_1^c,b_2^c,[b_1^c,b_2^c]}$ holds, we have $\C{D_{n_3+1}}=\Span{\partial_{x_3},b_1^c}\neq D_{n_3}$. For $n_1\geq 1$ or $[a_2,b_2^c]\neq 0\Mod\Span{\partial_{x_3},b_1^c,b_2^c,[b_1^c,b_2^c]}$, the distribution $D_{n_3+1}$ explicitly depends on $x_{3,1}^1$, \ie $[\partial_{x_{3,1}^1},D_{n_3+1}]\not\subset D_{n_3+1}$. In this case, because of $\partial_{x_{3,1}^1}\in D_{n_3}$, the condition $\C{D_{n_3+1}}\neq D_{n_3}$ also holds.
			\end{factProof}
			\vspace{-2em}
			\subsection*{Details sufficiency}
			 	\paragraph*{Details step 1.} Consider again the sequence of nested involutive distributions \eqref{eq:involutiveSequenceProof}. The dimensions of these distributions are as follows. By assumption, we have $\Dim{D_i}=2i$ for $i=1,\ldots,n_3+1$ and by construction, we have $\Dim{\Delta_0}=\Dim{D_{n_3-1}}+1=2n_3-1$. Furthermore, by construction we have $\Dim{\Delta_1}=\Dim{D_{n_3}}+1=2n_3+1=\Dim{\Delta_0}+2$. Because of item \ref{it:a} and \ref{it:b}, Lemma \ref{lem:cartan} applies to $\Delta_1$. Thus, for the Cauchy characteristic distributions $\C{\Delta_1^{(i)}}$, we have $ \Dim{\C{\Delta_1^{(i)}}}=\Dim{\Delta_0}+i$ for $i=1,\ldots,n_2-3$. For the involutive closure of $\Delta_1$, we have $\Dim{\overline{\Delta}_1}=\Dim{\Delta_0}+n_2$. Provided that an $x_1$-subsystem exists, \ie provided that $\overline{\Delta}_1\neq\mathcal{T}(\mathcal{X})$, we have either $\Dim{G_1}=\Dim{\overline{\Delta}_1}+2$ and $1\leq\Dim{G_{i+1}}-\Dim{G_i}\leq2$ or $\Dim{G_1}=\Dim{\overline{\Delta}_1}+1$ and $\Dim{G_{i+1}}=\Dim{G_i}+1$ for $i=1,\ldots,s-1$. In Step 1 of the sufficiency part of the proof of Theorem \ref{thm:1}, we claimed that after straightening out the distributions \eqref{eq:involutiveSequenceProof}, \ie applying a change of coordinates such that
				\begin{align*}
					\begin{aligned}
						D_1&=\Span{\partial_{x_3^{2n_3-1}},\partial_{x_3^{2n_3-2}}}\\
						&\vdotswithin{=}\\
						D_{n_3-1}&=\Span{\partial_{x_3^{2n_3-1}},\ldots,\partial_{x_3^2}}\\
						\Delta_0&=\Span{\partial_{x_3}}\\
						\C{\Delta_1^{(1)}}&=\Span{\partial_{x_3},\partial_{x_2^{n_2}}}\\
						&\vdotswithin{=}\\
						\C{\Delta_1^{(n_2-3)}}&=\Span{\partial_{x_3},\partial_{x_2^{n_2}},\ldots,\partial_{x_2^{4}}}\\
						\overline{\Delta}_1&=\Span{\partial_{x_3},\partial_{x_2}}\\
						G_1&=\Span{\partial_{x_3},\partial_{x_2},\partial_{x_1^{n_1}},\partial_{x_1^{n_1-1}}}\\
						&\vdotswithin{=}\\
						G_s&=\Span{\partial_{x_3},\partial_{x_2},\partial_{x_1}}=\mathcal{T}(\mathcal{X})\,,
					\end{aligned}
				\end{align*}
				the system is decomposed into the form
				\begin{align*}
					\begin{aligned}
						\dot{x}_1&=f_1(x_1,x_2^1,x_2^2,x_2^3)\\
						\dot{x}_2&=f_2(x_1,x_2,x_3^1,x_3^2,x_3^3)\\
						\dot{x}_3&=f_3(x_1,x_2,x_3,u^1,u^2)\,,
					\end{aligned}
				\end{align*}
				with $\Rank{\partial_{(x_2^1,x_2^2,x_2^3)}f_1}\leq2$, $\Rank{\partial_{(x_3^1,x_3^2,x_3^3)}f_2}=2$ and $\Rank{\partial_{(x_3^2,x_3^3)}f_2}=1$. In the following, we explain why this is indeed the case.
				\begin{propProof}[\ref{prop:1}]
					That $f_1$ only depends on the states $(x_1,x_2^1,x_2^2,x_2^3)$ and is independent of the states $(x_2^4,\ldots,x_2^{n_2},x_3)$ is implied by item \ref{it:c} condition \eqref{eq:compatibilityThm1} evaluated for $i=n_2-3$, \ie $[a,\C{\Delta_1^{(n_2-3)}}]\subset\Delta_1^{(n_2-3)}$. We have $\C{\Delta_1^{(n_2-3)}}=\Span{\partial_{x_3},\partial_{x_2^{n_2}},\ldots,\partial_{x_2^{4}}}$. If $f_1$ would depend on any of the states $(x_2^4,\ldots,x_2^{n_2},x_3)$, then $[a,\C{\Delta_1^{(n_2-3)}}]$ would contain vector fields with a $\partial_{x_1}$-component and thus, $[a,\C{\Delta_1^{(n_2-3)}}]$ would not be contained in $\Delta_1^{(n_2-3)}\subset\overline{\Delta}_1=\Span{\partial_{x_3},\partial_{x_2}}$.\\\\If $n_2=3$, condition \eqref{eq:compatibilityThm1} does not exist and we have $x_2=(x_2^1,x_2^2,x_2^3)$. In this case, $f_1$ can depend on all the states $x_2$. To show that also in this case $f_1$ cannot depend on the stats $x_3$, note that we have $\Delta_0=\Span{\partial_{x_3}}$ and by construction $[a,\Delta_0]\subset\Delta_1\subset\overline{\Delta}_1=\Span{\partial_{x_3},\partial_{x_2}}$. If $f_1$ would depend on $x_3$, then $[a,\Delta_0]$ would contain vector fields with a $\partial_{x_1}$-components and thus, $[a,\Delta_0]$ would not be contained in $\overline{\Delta}_1=\Span{\partial_{x_3},\partial_{x_2}}$.\\\\
					Next, let us show that $f_2=f_2(x_1,x_2,x_3^1,x_3^2,x_3^3)$. We have $D_{n_3-2}=\Span{\partial_{x_3^{n_3}},\ldots,\partial_{x_3^4}}$ and $D_{n_3-1}=\Span{\partial_{x_3^{n_3}},\ldots,\partial_{x_3^2}}$. Since by construction $[a,D_{n_3-2}]\subset D_{n_3-1}$ and $D_{n_3-1}$ contains no vector fields which have a $\partial_{x_2}$-component, it follows that $f_2$ is indeed independent of the states $(x_3^4,\ldots,x_3^{2n_3-1})$. 
				\end{propProof}
				\begin{propProof}[\ref{prop:2}]
					We have to show several rank conditions. First, let us show that $\Rank{\partial_{(x_2^1,x_2^2,x_2^3)}f_1}\leq 2$ holds, or equivalently that $\Dim{G_1}\leq\Dim{\overline{\Delta}_1}+2$ holds. Recall that because of item \ref{it:a} and \ref{it:b}, Lemma \ref{lem:cartan} applies to $\Delta_1$ and thus, we have $\Dim{\C{\Delta_1^{(n_2-3)}}}=\Dim{\overline{\Delta}_1}-3$ and $\Dim{\Delta_1^{(n_2-3)}}=\Dim{\overline{\Delta}_1}-1$. Therefore, there exist three vector fields $v_1,\ldots,v_3$ such that $\Delta_1^{(n_2-3)}=\C{\Delta_1^{(n_2-3)}}+\Span{v_1,v_2}$ and $\overline{\Delta}_1^{(n_2-3)}=\C{\Delta_1^{(n_2-3)}}+\Span{v_1,v_2,v_3}$. Due to item \ref{it:c} condition \eqref{eq:coupling}, \ie $\Dim{\overline{\Delta}_1+[a,\Delta_1^{(n_2-3)}]}=\Dim{\overline{\Delta}_1}+1$, the vector fields $[a,v_1]$ and $[a,v_2]$ are collinear$\Mod\overline{\Delta}_1$. Therefore, the dimension of $G_1=\overline{\Delta}_1+\Span{[a,v_1],[a,v_2],[a,v_3]}$ exceeds that of $\overline{\Delta}_1$ at most by two and thus $\Rank{\partial_{(x_2^1,x_2^2,x_2^3)}f_1}\leq2$ holds.\\\\
					Next, let us show that $\Rank{\partial_{(x_3^1,x_3^2,x_3^3)}f_2}=2$ holds. Note that because of $f_2=f_2(x_1,x_2,x_3^1,x_3^2,x_3^3)$ we have $\Rank{\partial_{(x_3^1,x_3^2,x_3^3)}f_2}=\Rank{\partial_{x_3}f_2}$, \ie since $f_2$ is independent of $(x_3^4,\ldots,x_3^{2n_3-1})$, the Jacobian matrices of $f_2$ with respect to $(x_3^1,x_3^2,x_3^3)$ and with respect to all the states $x_3$ have the same rank. We have $\Delta_0=\Span{\partial_{x_3}}$. From $\Delta_1=\Delta_0+[a,\Delta_0]$ and $\Dim{\Delta_1}=\Dim{\Delta_0}+2$ it follows that $\Rank{\partial_{x_3}f_2}=2$ holds and thus, also $\Rank{\partial_{(x_3^1,x_3^2,x_3^3)}f_2}=2$ indeed holds\footnote{The relation $\Delta_1=\Delta_0+[a,\Delta_0]$ actually only holds for $n_3\geq 2$, for $n_3=1$, we have $\Dim{\Delta_0}=1$ and thus, the distribution $\Delta_1=\Delta_0+[a,\Delta_0]$ would be of dimension 2. However, in this case, extending all the distributions by $\Span{\partial_{u^1},\partial_{u^2}}$ (which is the same as setting $D_{n_3-1}=D_0=\Span{\partial_{u^1},\partial_{u^2}}$), it can be shown that $f_2=f_2(x_1,x_2,x_3^1,u^1,u^2)$ and that the rank conditions $\Rank{\partial_{(x_3^1,u^1,u^2)}f_2}=2$ and $\Rank{\partial_{(u^1,u^2)}f_2}=1$ hold.}.\\\\					
					The last rank condition, namely $\Rank{\partial_{(x_3^2,x_3^3)}f_2}=1$, follows from $D_{n_3-1}\subset\Delta_0\subset D_{n_3}$, $\Dim{D_{n_3}}=\Dim{\Delta_0}+1$ and $D_{n_3}=D_{n_3-1}+[a,D_{n_3-1}]$. Since $D_{n_3-1}\subset\Delta_0\subset D_{n_3}$, we also have $D_{n_3}=\Delta_0+[a,D_{n_3-1}]$ and since $\Dim{D_{n_3}}=\Dim{\Delta_0}+1$, $[a,D_{n_3-1}]$ yields one direction which is not already contained in $\Delta_0=\Span{\partial_{x_3}}$. Because of $D_{n_3-1}=\Span{\partial_{x_3^{2n_3-1}},\ldots,\partial_{x_3^3},\partial_{x_3^2}}$, this implies $\Rank{\partial_{(x_3^2,x_3^3)}f_2}=1$.
				\end{propProof}
			\vspace{-1em}
			\paragraph*{Details step 2.} In the following we show the Propositions \ref{prop:3}, \ref{prop:4} and \ref{prop:5}, \ie we show that in any case regarding the actual form of the $x_1$-subsystem, there always exists a suitable flat output of the $x_2$-subsystem which is compatible with its (extended) chained structure. In Proposition \ref{prop:3}, the case $\Rank{\partial_{(x_2^1,x_2^2,x_2^3)}f_1}=2$ is addressed. In this case, the $x_1$-subsystem consists of two integrator chains and we have to show that the functions $\varphi^j(\bar{x}_1,x_2^1,x_2^2,x_2^3)$, $j=1,2$ in \eqref{eq:step2}, \ie the "inputs" of the integrator chains, form a compatible flat output of the $x_2$-subsystem. In Proposition \ref{prop:4}, the case $\Rank{\partial_{(x_2^1,x_2^2,x_2^3)}f_1}=1$ is addressed. In this case, the $x_1$-subsystem consists only of one integrator chain, it determines only one function $\varphi^1(\bar{x}_1,x_2^1,x_2^2,x_2^3)$ and we have to show that there always exists a second function $\varphi^2(\bar{x}_1,x_2^1,x_2^2,x_2^3)$, which together with the function $\varphi^1(\bar{x}_1,x_2^1,x_2^2,x_2^3)$, forms a compatible flat output of the $x_2$-subsystem. The case $\overline{\Delta}_1=\mathcal{T}(\mathcal{X})$, in which no $x_1$-subsystem exists, is addressed in Proposition \ref{prop:5}. In this case, we only have to construct a pair of functions which forms a compatible flat output of the $x_2$-subsystem, the flat output need not fulfill additional properties imposed by an $x_1$-subsystem. The results of \cite{LiRespondek:2012} regarding the normal chained form in fact directly apply. In all of these cases the existence of a certain involutive distribution $L$ is shown, the distribution $L$ is of importance in the problem of transforming the $x_2$-subsystem into (extended) chained form. 
			\begin{propProof}[\ref{prop:3}]
				For $\Rank{\partial_{(x_2^1,x_2^2,x_2^3)}f_1}=2$, the $x_1$-subsystem consists of two integrator chains and we have to show that the functions $\varphi^j(\bar{x}_1,x_2^1,x_2^2,x_2^3)$, $j=1,2$ in \eqref{eq:step2} meet $L=(\Span{\D\bar{x}_1,\D\varphi^1,\D\varphi^2})^\perp\subset\Delta_1^{(n_2-3)}$. To show this, let us first construct a special basis for the distribution $\Delta_1^{(n_2-3)}$. Recall that Lemma \ref{lem:cartan} applies to $\Delta_1$ and thus, we have $\Dim{\C{\Delta_1^{(n_2-3)}}}=\Dim{\Delta_1^{(n_2-3)}}-2$. Therefore, the distribution $\Delta_1^{(n_2-3)}$ can be represented as $\Delta_1^{(n_2-3)}=\C{\Delta_1^{(n_2-3)}}+\Span{v_1,v_2}$, with suitable vector fields $v_1,v_2\in\Delta^{(n_2-3)}_1$. Because of item \ref{it:c} condition \eqref{eq:coupling}, \ie $\Dim{\overline{\Delta}_1+[a,\Delta_1^{(n_2-3)}]}=\Dim{\overline{\Delta}_1}+1$, the vector fields $[a,v_1]$ and $[a,v_2]$ are collinear$\Mod\overline{\Delta}_1$, \ie (permute $v_1$ and $v_2$ if necessary) $[a,v_2]=\lambda[a,v_1]\Mod\overline{\Delta}_1$. The vector field $\tilde{v}_2=v_2-\lambda v_1$ therefore satisfies 
				\begin{align*}
					\begin{aligned}
						[a,\tilde{v}_2]&=\lambda[a,v_1]-\lambda[a,v_1]-\Lie_a\lambda v_1=0\Mod\overline{\Delta}_1
					\end{aligned}
				\end{align*}
				and thus, $\Lie_{\tilde{v}_2}\varphi^j=0$. Therefore, by choosing a basis of $\C{\Delta_1^{(n_2-3)}}$ together with $v_1$ and $\tilde{v}_2$ as basis for $\Delta_1^{(n_2-3)}$, \ie $\Delta_1^{(n_2-3)}=\C{\Delta_1^{(n_2-3)}}+\Span{v_1,\tilde{v}_2}$, we obtain a basis of which all basis vector fields, except for $v_1$, are annihilated by $\D\varphi^1$ and $\D\varphi^2$. The 1-form 
				\begin{align*}
					\begin{aligned}
						\omega&=(\D\varphi^2\rfloor v_1)\D\varphi^1-(\D\varphi^1\rfloor v_1)\D\varphi^2\,,
					\end{aligned}
				\end{align*}
				because of
				\begin{align*}
					\begin{aligned}
						\omega\rfloor v_1&=(\D\varphi^2\rfloor v_1)(\D\varphi^1\rfloor v_1)-(\D\varphi^1\rfloor v_1)(\D\varphi^2\rfloor v_1)=0\,,
					\end{aligned}
				\end{align*}
				annihilates all basis vector fields of $\Delta_1^{(n_2-3)}$. This 1-form together with $\D\bar{x}_1$ therefore spans the annihilator of $\Delta_1^{(n_2-3)}$. This shows that the annihilator of $\Delta_1^{(n_2-3)}$ is indeed a sub-codistribution of $L^\perp=\Span{\D\bar{x}_1,\D\varphi^1,\D\varphi^2}$, \ie $(\Delta_1^{(n_2-3)})^\perp\subset\Span{\D\bar{x}_1,\D\varphi^1,\D\varphi^2}=L^\perp$, or, equivalently $L=(\Span{\D\bar{x}_1,\D\varphi^1,\D\varphi^2})^\perp\subset\Delta_1^{(n_2-3)}$ indeed holds.
			\end{propProof}
			\begin{propProof}[\ref{prop:4}]
				For $\Rank{\partial_{(x_2^1,x_2^2,x_2^3)}f_1}=1$, the $x_1$-subsystem consists only of one integrator chain. In this case, the $x_1$-subsystem determines one function $\varphi^1(\bar{x}_1,x_2^1,x_2^2,x_2^3)$, \ie the "input" of the single integrator chain. We have to show that in this case, there always exists a second function $\varphi^2(\bar{x}_1,x_2^1,x_2^2,x_2^3)$ which together with the function $\varphi^1(\bar{x}_1,x_2^1,x_2^2,x_2^3)$, forms a compatible flat output of the $x_2$-subsystem. For that, we again make use of a special basis for the distribution $\Delta_1^{(n_2-3)}$, namely $\Delta_1^{(n_2-3)}=\C{\Delta_1^{(n_2-3)}}+\Span{v_1,v_2}$ with $v_j=v_j^i(x_2^1,x_2^2,x_2^3)\partial_{x_2^i}$, $i=1,\ldots,3$, and thus $[c,v_j]=0\Mod\C{\Delta_1^{(n_2-3)}}$ for any $c\in\C{\Delta_1^{(n_2-3)}}$. 
				The non-zero vector field $\tilde{v}=(\D\varphi^1\rfloor v_2)v_1-(\D\varphi^1\rfloor v_1)v_2$ annihilates $\D\varphi^1$, \ie $\Lie_{\tilde{v}}\varphi^1=0$. Together with a basis of  $\C{\Delta_1^{(n_2-3)}}$, this vector field $\tilde{v}$ spans the involutive distribution $L=\C{\Delta_1^{(n_2-3)}}+\Span{\tilde{v}}\subset\Delta_1^{(n_2-3)}$. We obviously have $\D\varphi^1\in L^\perp$ and since $L$ is involutive, there exists a second function $\varphi^2(x_2^1,x_2^2,x_2^3)$ such that $L=(\Span{\D\bar{x}_1,\D\varphi^1,\D\varphi^2})^\perp\subset\Delta_1^{(n_2-3)}$. The distribution $L$ is uniquely determined by the function $\varphi^1(x_2^1,x_2^2,x_2^3)$, the particular choice of the basis vector fields $v_1$ and $v_2$ does not matter. In fact, for any pair of vector fields $w_1,w_2$, which together with a basis of $\C{\Delta_1^{(n_2-3)}}$ spans $\Delta_1^{(n_2-3)}$, the distribution $\tilde{L}=\C{\Delta_1^{(n_2-3)}}+\Span{(\D\varphi^1\rfloor w_2)w_1-(\D\varphi^2\rfloor w_1)w_2}$ coincides with $L=\C{\Delta_1^{(n_2-3)}}+\Span{(\D\varphi^1\rfloor v_2)v_1-(\D\varphi^2\rfloor v_1)v_2}$ from above, \ie $\tilde{L}=L$. To show this, note that any vector fields $w_1$ and $w_2$, which together with $\C{\Delta_1^{(n_2-3)}}$ span the distribution $\Delta_1^{(n_2-3)}$, can be written as a linear combination $w_j=\beta_j^1v_1+\beta_j^2v_2\Mod\C{\Delta_1^{(n_2-3)}}$, $j=1,2$, with the vector fields $v_j=v_j^i(x_2^1,x_2^2,x_2^3)\partial_{x_2^i}$, $i=1,\ldots,3$ from above and functions $\beta_j^i=\beta_j^i(\bar{x}_1,x_2,x_3)$ and $\beta_1^1\beta_2^2-\beta_2^1\beta_1^2\neq 0$. Because of
				\begin{align*}
					\begin{aligned}
						(\D\varphi^1\rfloor w_2)w_1-(\D\varphi^2\rfloor w_1)w_2&=\left((\D\varphi^1\rfloor v_1)\beta_2^1+(\D\varphi^1\rfloor v_2)\beta_2^2\right)(\beta_1^1v_1+\beta_1^2v_2)-\\
						&\hspace{3em}\left((\D\varphi^1\rfloor v_1)\beta_1^1+(\D\varphi^1\rfloor v_2)\beta_1^2\right)(\beta_2^1v_1+\beta_2^2v_2)\Mod\C{\Delta_1^{(n_2-3)}}\\
						&=\underbrace{(\beta_1^1\beta_2^2-\beta_2^1\beta_1^2)}_{\neq 0}\left((\D\varphi^1\rfloor v_2)v_1-(\D\varphi^1\rfloor v_1)v_2\right)\Mod\C{\Delta_1^{(n_2-3)}}\,,
					\end{aligned}
				\end{align*}
				the distribution $\tilde{L}=\C{\Delta_1^{(n_2-3)}}+\Span{(\D\varphi^1\rfloor w_2)w_1-(\D\varphi^2\rfloor w_1)w_2}$ indeed coincides with $L=\C{\Delta_1^{(n_2-3)}}+\Span{(\D\varphi^1\rfloor v_2)v_1-(\D\varphi^2\rfloor v_1)v_2}$. Therefore, we obtain the unique distribution $L$ via $L=\C{\Delta_1^{(n_2-3)}}+\Span{(\D\varphi^1\rfloor w_2)w_1-(\D\varphi^2\rfloor w_1)w_2}$ by choosing arbitrary vector fields $w_1$ and $w_2$, which together with $\C{\Delta_1^{(n_2-3)}}$ span the distribution $\Delta_1^{(n_2-3)}$.
			\end{propProof}
			\begin{propProof}[\ref{prop:5}]
				For $\overline{\Delta}_1=\mathcal{T}(\mathcal{X})$, there does not exist an $x_1$-subsystem. We have to show that in this case, there always exist two functions $\varphi^j(x_2^1,x_2^2,x_2^3)$, $j=1,2$ which fulfill $L=(\Span{\D\varphi^1,\D\varphi^2})^\perp\subset\Delta_1^{(n_2-3)}$. The construction of $L$ is essentially the same as in the proof of Proposition \ref{prop:4}. We just have to omit $\bar{x}_1$ and since there is no $x_1$-subsystem which determines a function $\varphi^1$, we have to choose one. A valid choice is any function $\varphi^1(x_2^1,x_2^2,x_2^3)$, $\D\varphi^1\neq 0$. The distribution $L$ is then obtained via $L=\C{\Delta_1^{(n_2-3)}}+\Span{(\D\varphi^1\rfloor w_2)w_1-(\D\varphi^2\rfloor w_1)w_2}$, again with arbitrary vector fields $w_1$ and $w_2$, which together with $\C{\Delta_1^{(n_2-3)}}$ span the distribution $\Delta_1^{(n_2-3)}$. The distribution $L$ is again independent of the particular choice of $w_1$ and $w_2$, it is uniquely determined by the choice of $\varphi^1(x_2^1,x_2^2,x_2^3)$. This construction in fact coincides with the construction of $L$ provided in \cite{LiRespondek:2012}, Theorem 2.10, \ie in case that no $x_1$-subsystem exists, the results from \cite{LiRespondek:2012} directly apply.
			\end{propProof}
			There is another way to calculate the distribution $L$ to a given function $\varphi^1$ (either determined by the $x_1$-subsystem or chosen if no $x_1$-subsystem exists). As we have seen above, the distribution $L$ is uniquely determined by the function $\varphi^1$ and there always exists a second function $\varphi^2$ such that $L=(\Span{\D\bar{x}_1,\D\varphi^1,\D\varphi^2})^\perp\subset\Delta_1^{(n_2-3)}$ (omit $\bar{x}_1$ in case that there is no $x_1$-subsystem). Recall that we have $\Dim{\Delta_1^{(n_2-3)}}=\Dim{\overline{\Delta}_1}-1$ and $\overline{\Delta}_1=\Span{\partial_{x_3},\partial_{x_2}}$. The annihilator of $\Delta_1^{(n_2-3)}$ is thus of the form $(\Delta_1^{(n_2-3)})^\perp=\Span{\D\bar{x}_1,\omega}$ and because of $L=(\Span{\D\bar{x}_1,\D\varphi^1,\D\varphi^2})^\perp\subset\Delta_1^{(n_2-3)}$, the annihilator of $\Delta_1^{(n_2-3)}$ is a sub-codistribution of $L^\perp=\Span{\D\bar{x}_1,\D\varphi^1,\D\varphi^2}$, \ie $(\Delta_1^{(n_2-3)})^\perp=\Span{\D\bar{x}_1,\omega}\subset\Span{\D\bar{x}_1,\D\varphi^1,\D\varphi^2}$. Thus, the 1-form $\omega$ is a linear combination of the differentials $\D\varphi^1$ and $\D\varphi^2$ and thus, we have $L^\perp=\Span{\D\bar{x}_1,\D\varphi^1,\omega}$. Therefore, given $\varphi^1$, we immediately obtain the annihilator of the associated distribution $L$ via $L^\perp=(\Delta_1^{(n_2-3)})^\perp+\Span{\D\varphi^1}$. By integrating this codistribution, we obtain a possible second function $\varphi^2$.
			\paragraph*{Details step 4.} In the following, we show that after applying the transformation \eqref{eq:firstInput}, the $x_2$-subsystem takes the form \eqref{eq:subsys2triangular}, as asserted in Proposition \ref{prop:6}. The Fact \ref{ft:2} is shown subsequently.
			\begin{propProof}[\ref{prop:6}]
				The transformation \eqref{eq:firstInput} normalizes the first equation of the $x_2$-subsystem, \ie applying this transformation immediately yields an $x_2$-subsystem of the form
				\begin{align}\label{eq:step4_ap1}
					f_2:\quad\begin{aligned}
						\dot{\tilde{x}}_2^1&=x_{3,2}^1\\
						\dot{\tilde{x}}_2^2&=\tilde{f}_2^2(\bar{x}_1,\tilde{x}_2,x_3^1,x_{3,2}^1,x_3^3)\\
						&\vdotswithin{=}\\
						\dot{x}_2^{n_2}&=\tilde{f}_2^{n_2}(\bar{x}_1,\tilde{x}_2,x_3^1,x_{3,2}^1,x_3^3)\,.
					\end{aligned}
				\end{align}
				The same way as in the proof of Proposition \ref{prop:2} the rank condition $\Rank{\partial_{(x_3^2,x_3^3)}f_2}=1$ is shown, for \eqref{eq:step4_ap1} the rank condition $\Rank{\partial_{(x_{3,2}^1,x_3^3)}\tilde{f}_2}=1$ can be shown, which implies that the functions $\tilde{f}_2^j$, $j=2,\ldots,n_2$ are actually independent of $x_3^3$. In other words, the transformation \eqref{eq:firstInput} eliminates the redundancy among the inputs of the $x_2$-subsystem. For $n_3\geq 2$, we have $\Delta_1=\Delta_0+[a,\Delta_0]$ and thus $\Delta_1=\Span{\partial_{x_3},v_1,v_2}$ with the vector fields $v_1=\partial_{\tilde{x}_2^1}+\partial_{x_{3,2}^1}\tilde{f}_2^2\partial_{\tilde{x}_2^2}+\ldots+\partial_{x_{3,2}^1}\tilde{f}_2^{n_2}\partial_{x_2^{n_2}}$, $v_2=\partial_{x_3^1}\tilde{f}_2^2\partial_{\tilde{x}_2^2}+\ldots+\partial_{x_3^1}\tilde{f}_2^{n_2}\partial_{x_2^{n_2}}$. Furthermore, Lemma \ref{lem:cartan} applies to $\Delta_1$. Therefore, there exists a linear combination of $v_1$ and $v_2$ which is contained in $\C{\Delta_1^{(1)}}=\Span{\partial_{x_3},\partial_{x_2^{n_2}}}$\footnote{For $n_2=3$, we would have $\Delta_1^{(1)}=\overline{\Delta}_1$ and thus $\C{\Delta_1^{(1)}}=\overline{\Delta}_1$. In this case, replace $\C{\Delta_1^{(1)}}$ by $L=\Span{\partial_{x_3},\partial_{x_2^{n_2}}}$. Because of $L\subset\Delta_1$, there again exists a linear combination of $v_1$ and $v_2$ which is contained in $\Span{\partial_{x_3},\partial_{x_2^{n_2}}}$, where $\partial_{x_2^{n_2}}=\partial_{x_2^3}$ in this case.}. The vector field $v_1$ has a non-zero component in the $\partial_{\tilde{x}_2^1}$-direction, the vector field $v_2$ has not. Therefore, linear combinations of $v_1$ and $v_2$ which are contained in $\Span{\partial_{x_3},\partial_{x_2^{n_2}}}$ consist of $v_2$ only, and thus, we have $v_2\in\C{\Delta_1^{(1)}}$. Since $v_2=\partial_{x_3^1}\tilde{f}_2^2\partial_{\tilde{x}_2^2}+\ldots+\partial_{x_3^1}\tilde{f}_2^{n_2}\partial_{x_2^{n_2}}$, it follows that $x_3^1$ can only occur in the function $\tilde{f}_2^{n_2}$. Furthermore, in order for $\C{\Delta_1}=\Delta_0$ to hold, $x_{3,2}^1$ must occur affine in the functions $\tilde{f}_2^j$, $j=2,\ldots,n_2-1$, \ie the $x_2$-subsystem is actually of the form\footnote{If $n_3=1$, the state transformation \eqref{eq:firstInput} is replaced by the input transformation $\tilde{u}^2=\tilde{f}_2^1(\bar{x}^1,\tilde{x}_2,x_3^1,u^1,u^2)$, see footnote \ref{fn:transformation}. In \eqref{eq:step4_ap1}, $x_{3,2}^1$ and $x_3^3$ would then be replaced by $\tilde{u}^2$ and $u^1$. By an analogous reasoning as above, we would then find that the $x_2$-subsystem is actually independent of $u^1$ and that $x_3^3$ again occurs only in the very last equation of the $x_2$-subsystem, \ie the $x_2$-subsystem would again be of the form \eqref{eq:step4_ap2}, with $x_{3,2}^1$ replaced by $\tilde{u}^2$. That the input $\tilde{u}^2$ occurs affine in the $x_2$-subsystem would follow directly from the fact that we started with an AI-system and only applied transformations which preserve the AI structure.}
				\begin{align}\label{eq:step4_ap2}
					f_2:\quad\begin{aligned}
						\dot{\tilde{x}}_2^1&=x_{3,2}^1\\
						\dot{\tilde{x}}_2^2&=b_2^2(\bar{x}_1,\tilde{x}_2)x_{3,2}^1+a_2^2(\bar{x}_1,\tilde{x}_2)\\
						&\vdotswithin{=}\\
						\dot{x}_2^{n_2-1}&=b_2^{n_2-1}(\bar{x}_1,\tilde{x}_2)x_{3,2}^1+a_2^{n_2-1}(\bar{x}_1,\tilde{x}_2)\\
						\dot{x}_2^{n_2}&=g(\bar{x}_1,\tilde{x}_2,x_3^1,x_{3,2}^1)\,.
					\end{aligned}
				\end{align}
				Next, we show that the functions $b_2^i$ in \eqref{eq:step4_ap2} depend on the states of the $x_2$-subsystem in a triangular manner. Lemma \ref{lem:cartan} applies to $\Delta_1$. Based on that, we will first show that $\Delta_1^{(i)}=\Span{b_2}+\C{\Delta_1^{(i+1)}}$, $i=0,\ldots,n_2-3$, with $b_2=\partial_{\tilde{x}_2^1}+b_2^2\partial_{\tilde{x}_2^2}+b_2^3\partial_{x_2^3}+\ldots+b_2^{n_2-1}\partial_{x_2^{n_2-1}}$, \ie that the derived flags $\Delta_1^{(i)}$ are composed of the one-dimensional distribution spanned by the vector field $b_2$ and the Cauchy characteristic distributions of their next derived flags. To show this, note that the vector field $b_2$ has a component in the $\partial_{\tilde{x}_2^1}$-direction. Thus, it cannot belong to any of the Cauchy characteristics $\C{\Delta_1^{(i+1)}}=\Span{\partial_{x_3},\partial_{x_2^{n_2}},\ldots,\partial_{x_2^{n_2-i}}}$, $i=0,\ldots,n_2-4$. However, because of $b_2\in\Delta_1$, the vector field $b_2$ also belongs to all the derived flags of $\Delta_1$. Furthermore, because of Lemma \ref{lem:cartan}, we have $\Dim{\C{\Delta_1^{(i+1)}}}=\Dim{\Delta_1^{(i)}}-1$, $i=0,\ldots,n_2-4$. Thus, $\Span{b_2}$ completes $\C{\Delta_1^{(i+1)}}$ to $\Delta_1^{(i)}$. By construction, we furthermore have $L=(\Span{\D\bar{x}_1,\D\varphi^1,\D\varphi^2})^\perp\subset\Delta_1^{(n_2-3)}$ and $b_2\notin L$\footnote{Note that we have $L=\Span{\partial_{x_3},\partial_{x_2^{n_2}},\ldots,\partial_{x_2^3}}$. Thus, the vector field $b_2$, by having a component in the $\partial_{\bar{x}_2^1}$-direction cannot be contained in $L$.}. Thus, $\Delta_1^{(n_2-3)}=\Span{b_2}+L$ holds, \ie $\Delta_1^{(n_2-3)}=\Span{\partial_{x_3},\partial_{x_2^{n_2}},\ldots,\partial_{x_2^3},b_2}$. In conclusion, we have $\Delta_1^{(i)}=\Span{\partial_{x_3},\partial_{x_2^{n_2}},\ldots,\partial_{x_2^{n_2-i}},b_2}$, $i=0,\ldots,n_2-3$, from which $\partial_{x_2^j}b_2^k=0$ for $k+2\leq j\leq n_2$, $k=2,\ldots,n_2-2$, and $\partial_{x_2^j}b_2^k\neq 0$ for $j=k+1$, $k=2,\ldots,n_2-1$ follows. This exactly describes the triangular dependence of the functions $b_2^i$, $i=2,\ldots,n_2-1$ on the states $(x_2^3,\ldots,x_2^{n_2})$ in \eqref{eq:subsys2triangular}.\\\\
				The triangular dependence of the functions $a_2^i$, $i=2,\ldots,n_2-1$ on the states $(x_2^3,\ldots,x_2^{n_2})$ in \eqref{eq:subsys2triangular} is implied by item \ref{it:c} condition \eqref{eq:compatibilityThm1}, \ie $[a,\C{\Delta_1^{(i)}}]\subset\Delta_1^{(i)}$, $i=1,\ldots,n_2-3$. We have $\C{\Delta_1^{(i)}}=\Span{\partial_{x_3},\partial_{x_2^{n_2}},\ldots,\partial_{x_2^{n_2-i+1}}}$, $i=1,\ldots,n_2-3$ and $\Delta_1^{(i)}=\Span{\partial_{x_3},\partial_{x_2^{n_2}},\ldots,\partial_{x_2^{n_2-i}},b_2}$, $i=0,\ldots,n_2-3$. Evaluating $[a,\C{\Delta_1^{(i)}}]\subset\Delta_1^{(i)}$, $i=1,\ldots,n_2-3$ therefore yields $\partial_{x_2^j}a_2^k=0$ for $k+2\leq j\leq n_2$, $k=2,\ldots,n_2-2$. The condition $[a,\C{\Delta_1^{(i)}}]\subset\Delta_1^{(i)}$, $i=1,\ldots,n_2-3$ in fact coincides with the compatibility condition \eqref{eq:compatibility} in Theorem \ref{thm:extendedChainedform} for the extended chained form.
			\end{propProof}
			\begin{factProof}[\ref{ft:2}]
				It follows from the construction of the $x_2$-subsystem that every component of the right hand side of the $x_2$-subsystem, \ie every function $\tilde{f}_2^i$, $i=1,\ldots,n_2$ in $\dot{\tilde{x}}_2=\tilde{f}_2(\bar{x}_1,\tilde{x}_2,x_3^1,x_3^2,x_3^3)$ explicitly depends on at least one of the inputs $(x_3^1,x_3^2,x_3^3)$ of the $x_2$-subsystem. Under the assumption that $\tilde{f}_2^1$ indeed explicitly depends on $x_3^2$, in Step 4 of the proof, $x_3^2$ is replaced by $x_{3,2}^1=\tilde{f}_2^1(\bar{x}_1,\tilde{x}_2,x_3^1,x_3^2,x_3^3)$, which results in an $x_2$-subsystem of the form \eqref{eq:subsys2triangular} and $D_{n_3-1}=\Span{\partial_{x_3^{2n_3-1}},\ldots,\partial_{x_3^3},\partial_{x_{3,2}^1}}$, \ie this transformation certainly keeps the distribution $D_{n_3-1}$ straightened out. To show that $\tilde{f}_2^1$ indeed depends on $x_3^2$ or $x_3^3$, let us instead replace the state $x_3^1$ by $x_{3,2}^1=\tilde{f}_2^1$ and keep $x_3^2$ as coordinate. By a similar reasoning as in the proof of Proposition \ref{prop:6}, it then follows that after this transformation, the $x_2$-subsystem takes the form
				\begin{align*}
					f_2:\quad\begin{aligned}
						\dot{\tilde{x}}_2^1&=x_{3,2}^1\\
						\dot{\tilde{x}}_2^2&=b_2^2(\bar{x}_1,\tilde{x}_2^1,\tilde{x}_2^2,x_2^3)x_{3,2}^1+a_2^2(\bar{x}_1,\tilde{x}_2^1,\tilde{x}_2^2,x_2^3)\\
						\dot{x}_2^3&=b_2^3(\bar{x}_1,\tilde{x}_2^1,\tilde{x}_2^2,x_2^3,x_2^4)x_{3,2}^1+a_2^3(\bar{x}_1,\tilde{x}_2^1,\tilde{x}_2^2,x_2^3,x_2^4)\\
						&\vdotswithin{=}\\
						\dot{x}_2^{n_2-1}&=b_2^{n_2-1}(\bar{x}_1,\tilde{x}_2)x_{3,2}^1+a_2^{n_2-1}(\bar{x}_1,\tilde{x}_2)\\
						\dot{x}_2^{n_2}&=g(\bar{x}_1,\tilde{x}_2,x_3^2,x_3^3,x_{3,2}^1)\,.
					\end{aligned}
				\end{align*}
				Furthermore, in the new coordinates, the distribution $D_{n_3-1}$ takes the form $D_{n_3-1}=\Span{\partial_{x_3^{2n_3-1}},\ldots,\partial_{x_3^4},\partial_{x_3^3}+\partial_{x_3^3}\tilde{f}_2^1\partial_{x_{3,2}^1},\partial_{x_3^2}+\partial_{x_3^2}\tilde{f}_2^1\partial_{x_{3,2}^1}}$. Assume that $\tilde{f}_2^1$ is independent of $x_3^2$ and $x_3^3$. Then, we have $D_{n_3-1}=\Span{\partial_{x_3^{2n_3-1}},\ldots,\partial_{x_3^4},\partial_{x_3^3},\partial_{x_3^2}}$, \ie $D_{n_3-1}$ is still straightened out. However, this leads to $D_{n_3}=D_{n_3-1}+[a,D_{n_3-1}]=\Span{\partial_{x_3},\partial_{x_2^{n_2}}}$ and $D_{n_3+1}=D_{n_3}+[a,D_{n_3}]=\Span{\partial_{x_3},\partial_{x_2^{n_2}},\partial_{x_2^{n_2-1}},\partial_{\tilde{x}_2^1}+b_2^2\partial_{\tilde{x}_2^2}+\ldots+b_2^{n_2-2}\partial_{x_2^{n_2-2}}}$ and in turn, $\C{D_{n_3+1}}=D_{n_3}$ would hold, or, if $n_2=3$, $D_{n_3+1}$ would be involutive. (In fact, this would lead exactly to the case mentioned in footnote \ref{fn:wrongInput}, where the longer integrator chain of the $x_3$-subsystem is attached to the "wrong" input of the $x_2$-subsystem.)\\\\
				Similarly, for $n_3=1$, it can be shown that an $\tilde{f}_2^1$ which does not depend on an input $u^1$ or $u^2$ leads to the same contradictions (see also footnote \ref{fn:transformation}). Assume that $\tilde{f}_2^1$ is independent of the inputs $u^1$ and $u^2$, \ie $\tilde{f}_2^1=\tilde{f}_2^1(\bar{x}_1,\tilde{x}_2,x_3^1)$. By applying the state transformation $\tilde{x}_3^1=\tilde{f}_2^1(\bar{x}_1,\tilde{x}_2,x_3^1)$, followed by a suitable input transformation, we would then obtain
				\begin{align*}
					f_2:&\quad\begin{aligned}
						\dot{\tilde{x}}_2^1&=\tilde{x}_3^1\\
						\dot{\tilde{x}}_2^2&=b_2^2(\bar{x}_1,\tilde{x}_2^1,\tilde{x}_2^2,x_2^3)\tilde{x}_3^1+a_2^2(\bar{x}_1,\tilde{x}_2^1,\tilde{x}_2^2,x_2^3)\\
						\dot{x}_2^3&=b_2^3(\bar{x}_1,\tilde{x}_2^1,\tilde{x}_2^2,x_2^3,x_2^4)\tilde{x}_3^1+a_2^3(\bar{x}_1,\tilde{x}_2^1,\tilde{x}_2^2,x_2^3,x_2^4)\\
						&\vdotswithin{=}\\
						\dot{x}_2^{n_2-1}&=b_2^{n_2-1}(\bar{x}_1,\tilde{x}_2)\tilde{x}_3^1+a_2^{n_2-1}(\bar{x}_1,\tilde{x}_2)\\
						\dot{x}_2^{n_2}&=\tilde{u}^2
					\end{aligned}\\
					f_3:&\quad\quad~\,\begin{aligned}
					\dot{\tilde{x}}_3^1&=\tilde{u}^1\,.
					\end{aligned}
				\end{align*}
				However, this would again lead to $D_{n_3}=D_1=\Span{\partial_{x_3},\partial_{x_2^{n_2}}}$ and $D_{n_3+1}=D_2=\Span{\partial_{x_3},\partial_{x_2^{n_2}},\partial_{x_2^{n_2-1}},\partial_{\tilde{x}_2^1}+b_2^2\partial_{\tilde{x}_2^2}+\ldots+b_2^{n_2-2}\partial_{x_2^{n_2-2}}}$ and therefore again to $\C{D_{n_3+1}}=D_{n_3}$ or an involutive $D_{n_3+1}$ in case that $n_2=3$.
			\end{factProof}
			\vspace{-2em}
		\subsection{Proof of the simple method for determining $b_p$}\label{ap:proofSimplifiedMethod}
			In the following we show why in the case $\Ad_a^{n_3+1}b_1\notin H$ or $\Ad_a^{n_3+1}b_2\notin H$, with the distribution $H=D_{n_3+1}+[D_{n_3},D_{n_3+1}]$, a vector field $b_p$ can indeed be determined from the criterion $\Ad_a^{n_3+1}b_p\in H$, as proposed in Remark \ref{rem:simplerMethod}. For a system of the form \eqref{eq:triangularFormDifferentLengths}, we obtain (recall that we have $D_{n_3}=\Span{\partial_{x_3},b_2^c+gb_1^c}$, $D_{n_3+1}=\Span{\partial_{x_3},b_1^c,b_2^c,[a,b_2^c+gb_1^c]}$ and $[a,b_2^c+gb_1^c]=[f_1,b_2^c]+\tilde{h}(x_1,x_2,x_{3,1}^1)[b_1^c,b_2^c]+[a_2,b_2^c]\Mod\Span{\partial_{x_3},b_1^c,b_2^c}$, see also proof of Fact \ref{ft:1})
			\begin{align}\label{eq:H}
				\begin{aligned}
					H&=D_{n_3+1}+[D_{n_3},D_{n_3+1}]\\
					&=\Span{\partial_{x_3},b_1^c,b_2^c,[a,b_2^c+gb_1^c],[b_2^c+gb_1^c,[a,b_2^c+gb_1^c]],[b_1^c,b_2^c]}\,.
				\end{aligned}
			\end{align}
			The input vector fields of \eqref{eq:triangularFormDifferentLengths} are $b_1=\partial_{x_{3,1}^{n_3}}$ and $b_2=\partial_{x_{3,2}^{n_3-1}}$, where $b_1$ is the input vector field belonging to the longer integrator chain of the $x_3$-subsystem. The conditions of Theorem \ref{thm:1} are met with any non-zero vector field $b_p$ which is collinear with $b_1$. For any vector field $b_p$ which is collinear with $b_1=\partial_{x_{3,1}^{n_3}}$, \ie  $b_p=\lambda\partial_{x_{3,1}^{n_3}}$ with an arbitrary non-zero function $\lambda$ of the state of the system, we obtain $\Ad_a^{n_3+1}b_p=\lambda\,(-1)^{n_3}\,[a,b_1^c]\Mod D_{n_3+1}$ and thus, because of $[a,b_1^c]\in\Span{\partial_{x_3},b_1^c,[b_1^c,b_2^c]}\subset H$ and $D_{n_3+1}\subset H$, we have $\Ad_a^{n_3+1}b_p\in H$. Whereas $\Ad_a^{n_3+1}b_2=(-1)^{n_3-1}[a,[a,b_2^c+gb_1^c]]$ may or may not be contained in $H$. If $\Ad_a^{n_3+1}b_2\notin H$, then $\Ad_a^{n_3+1}b_p\in H$ is indeed only met for vector fields $b_p$ which are collinear with $b_1$, \ie collinear with the input vector field of the longer integrator chain in the $x_3$-subsystem. However, if also $\Ad_a^{n_3+1}b_2\in H$, this criterion for determining a vector field $b_p$ is not applicable, since $\Ad_a^{n_3+1}b_p\in H$ would be met for every linear combination $b_p$ of the input vector fields of the system.
		\subsection{Analysis of the necessary condition (\ref{eq:alg})}\label{ap:analysisOfTheEquationsystem}
			In the following, we analyze the necessary condition \eqref{eq:alg} in terms of uniqueness of the direction of candidates for $b_p=\alpha^1b_1+\alpha^2b_2$. For a system of the form \eqref{eq:triangularFormDifferentLengths}, we have $b_1=\partial_{x_{3,1}^{n_3}}$ and $b_2=\partial_{x_{3,2}^{n_3-1}}$ and thus $v_1=\Ad_a^{n_3-1}b_1=(-1)^{n_3-1}\,\partial_{x_{3,1}^1}$ and $v_2=\Ad_a^{n_3-1}b_2=(-1)^{n_3-1}(b_2^c+gb_1^c)$, again with $b_1^c=\partial_{x_2^{n_2}}$ and $b_2^c=\partial_{x_2^1}+x_2^3\partial_{x_2^2}+\ldots+x_2^{n_2}\partial_{x_2^{n_2-1}}$. The conditions of Theorem \ref{thm:1} are met with any non-zero vector field $b_p$ which is collinear with $b_1$. Assume we apply a regular input transformation on the system. Then, we have $\tilde{b}_1=\beta_1^1b_1+\beta_1^2b_2$ and $\tilde{b}_2=\beta_2^1b_1+\beta_2^2b_2$, with $\beta_1^1\beta_2^2-\beta_1^2\beta_2^1\neq 0$, and $\tilde{a}=a+\gamma^1b_1+\gamma^2b_2$ and accordingly $\tilde{v}_1=\beta_1^1v_1+\beta_1^2v_2\Mod D_{n_3-1}$ and $\tilde{v}_2=\beta_2^1v_1+\beta_2^2v_2\Mod D_{n_3-1}$ ($\beta_i^j$ and $\gamma^j$ being functions of the state $x$ of the system, \ie $\beta_i^j=\beta_i^j(x)$ and $\gamma^j=\gamma^j(x)$). In the following, we show that by solving the necessary condition \eqref{eq:alg}, we obtain at most two non-collinear candidates for the vector field $b_p$, and that one of theses candidates is collinear with $b_1=\partial_{x_{3,1}^{n_3}}$. We start by inserting $\tilde{v}_1$ and $\tilde{v}_2$ into \eqref{eq:alg}, \ie
			\begin{align}\label{eq:alg_0}
				\begin{aligned}
					(\alpha^1)^2[\tilde{v}_1,[\tilde{a},\tilde{v}_1]]+2\alpha^1\alpha^2[\tilde{v}_1,[\tilde{a},\tilde{v}_2]]+(\alpha^2)^2[\tilde{v}_2,[\tilde{a},\tilde{v}_2]]&\overset{!}{\in}D_{n_3+1}\,.
				\end{aligned}
			\end{align}
			By inserting the corresponding expressions for $\tilde{v}_1$, $\tilde{v}_2$ and $\tilde{a}$ from above, we obtain\footnote{Note that because of $b_j,\tilde{v}_j\in D_{n_3}$ and the involutivity of $D_{n_3}$, we have $[\tilde{a},\tilde{v}_j]=[a+\gamma^1b_1+\gamma^2b_2,\tilde{v}_j]=[a,\tilde{v}_j]\Mod D_{n_3}$.}
			\begin{align*}
				\begin{aligned}
					&(\alpha^1)^2[\beta_1^1v_1+\beta_1^2v_2,[a,\beta_1^1v_1+\beta_1^2v_2]]+2\alpha^1\alpha^2[\beta_1^1v_1+\beta_1^2v_2,[a,\beta_2^1v_1+\beta_2^2v_2]]+\\
					&\hspace{19em}(\alpha^2)^2[\beta_2^1v_1+\beta_2^2v_2,[a,\beta_2^1v_1+\beta_2^2v_2]]\overset{!}{\in}D_{n_3+1}\,.
				\end{aligned}
			\end{align*}
			Expanding yields
			\begin{align*}
				\begin{aligned}
					&(\alpha^1)^2\left((\beta_1^1)^2[v_1,[a,v_1]]+\beta_1^1\beta_1^2[v_1,[a,v_2]]+\beta_1^2\beta_1^1[v_2,[a,v_1]]+(\beta_1^2)^2[v_2,[a,v_2]]\right)+\\
					&\hspace{2em}2\alpha^1\alpha^2\left(\beta_1^1\beta_2^1[v_1,[a,v_1]]+\beta_1^1\beta_2^2[v_1,[a,v_2]]+\beta_1^2\beta_2^1[v_2,[a,v_1]]+\beta_1^2\beta_2^2[v_2,[a,v_2]]\right)+\\
					&\hspace{2em}(\alpha^2)^2\left((\beta_2^1)^2[v_1,[a,v_1]]+\beta_2^1\beta_2^2[v_1,[a,v_2]]+\beta_2^2\beta_2^1[v_2,[a,v_1]]+(\beta_2^2)^2[v_2,[a,v_2]]\right)\overset{!}{\in}D_{n_3+1}\,.
				\end{aligned}	
			\end{align*}
			With $[v_2,[a,v_1]]=[v_1,[a,v_2]]\Mod D_{n_3+1}$ (following from the Jacobi identity), and $[v_1,[a,v_1]]\in D_{n_3+1}$ (actually $[v_1,[a,v_1]]=0$), we obtain
			\begin{align*}
				\begin{aligned}
					&(\alpha^1)^2\left(2\beta_1^1\beta_1^2[v_1,[a,v_2]]+(\beta_1^2)^2[v_2,[a,v_2]]\right)+\\
					&\hspace{2em}2\alpha^1\alpha^2\left((\beta_1^1\beta_2^2+\beta_1^2\beta_2^1)[v_1,[a,v_2]]+\beta_1^2\beta_2^2[v_2,[a,v_2]]\right)+\\
					&\hspace{2em}(\alpha^2)^2\left(2\beta_2^1\beta_2^2[v_1,[a,v_2]]+(\beta_2^2)^2[v_2,[a,v_2]]\right)\overset{!}{\in}D_{n_3+1}\,,
				\end{aligned}	
			\end{align*}
			and after some rearranging 
			\begin{align*}
				\begin{aligned}
					&2\left((\alpha^1)^2\beta_1^1\beta_1^2+\alpha^1\alpha^2(\beta_1^1\beta_2^2+\beta_1^2\beta_2^1)+(\alpha^2)^2\beta_2^1\beta_2^2\right)[v_1,[a,v_2]]+\\
					&\hspace{10em}\left((\alpha^1)^2(\beta_1^2)^2+2\alpha^1\alpha^2\beta_1^2\beta_2^2+(\alpha^2)^2(\beta_2^2)^2\right)[v_2,[a,v_2]]\overset{!}{\in}D_{n_3+1}\,,
				\end{aligned}	
			\end{align*}
			and finally
			\begin{align}\label{eq:alg_1}
				\begin{aligned}
					&(\alpha^1\beta_1^2+\alpha^2\beta_2^2)\left(2(\alpha^1\beta_1^1+\alpha^2\beta_2^1)[v_1,[a,v_2]]+(\alpha^1\beta_1^2+\alpha^2\beta_2^2)[v_2,[a,v_2]]\right)\overset{!}{\in}D_{n_3+1}\,.
				\end{aligned}	
			\end{align}
			In the following, we have to distinguish between two cases, namely between $[v_1,[a,v_2]]$ and $[v_2,[a,v_2]]$ being collinear$\Mod D_{n_3+1}$ or not.
			\paragraph*{Case 1:} Let us first consider the case $[v_1,[a,v_2]]$ and $[v_2,[a,v_2]]$ not being collinear$\Mod D_{n_3+1}$. In this case, there does not exist a non-trivial linear combination of the vector fields $[v_1,[a,v_2]]$ and $[v_2,[a,v_2]]$ which is contained in $D_{n_3+1}$. Furthermore, the factors $\alpha^1\beta_1^1+\alpha^2\beta_2^1$ and $\alpha^1\beta_1^2+\alpha^2\beta_2^2$ cannot vanish simultaneously for $\alpha^1\neq 0$ or $\alpha^2\neq 0$\footnote{All non-trivial solutions of the linear homogeneous equation $\alpha^1\beta_1^2+\alpha^2\beta_2^2$ are of the form $\alpha^1=\lambda\beta_2^2$, $\alpha^2=-\lambda\beta_1^2$ with arbitrary $\lambda\neq 0$. We have at least $\beta_2^2\neq 0$ or $\beta_1^2\neq 0$, otherwise, the input transformation from above would not be invertible (\ie for $\beta_1^1\beta_2^2-\beta_1^2\beta_2^1=0$, the new input vector fields $\bar{b}_1$ and $\bar{b}_2$ would be linearly dependent). Inserting this solution into the second factor $\alpha^1\beta_1^1+\alpha^2\beta_2^1$ yields  $\lambda(\beta_1^1\beta_2^2-\beta_1^2\beta_2^1)$. This term can only vanish for $\lambda=0$ since $\beta_1^1\beta_2^2-\beta_1^2\beta_2^1\neq 0$ for a regular transformation. However, $\lambda=0$ is the trivial solution $\alpha^1=\alpha^2=0$.}, or in other words, we cannot chose $\alpha^1$ and $\alpha^2$ such that in \eqref{eq:alg_1} there occurs a trivial linear combination of $[v_1,[a,v_2]]$ and $[v_2,[a,v_2]]$. Thus, in this case, in order for \eqref{eq:alg_1} to hold, the factor $(\alpha^1\beta_1^2+\alpha^2\beta_2^2)$ must vanish and thus $\alpha^1=\lambda\beta_2^2$ and $\alpha^2=-\lambda\beta_1^2$ with arbitrary $\lambda\neq 0$, \ie in this case the solution of the necessary condition \eqref{eq:alg_0} is unique up to a multiplication with arbitrary $\lambda\neq 0$. With this solution, for $b_p$ we obtain
			\begin{align*}
				\begin{aligned}
					b_p&=\alpha^1\tilde{b}_1+\alpha^2\tilde{b}_2\\
					&=\lambda(\beta_2^2(\beta_1^1b_1+\beta_1^2b_2)-\beta_1^2(\beta_2^1b_1+\beta_2^2b_2))\\
					&=\underbrace{\lambda(\beta_1^1\beta_2^2-\beta_1^2\beta_2^1)}_{\neq 0}b_1\,,
				\end{aligned}
			\end{align*}
			\ie we indeed recover the direction of $b_1=\partial_{x_{3,1}^{n_3}}$.
			\paragraph*{Case 2:} In this case, there exists a non-trivial linear combination of the vector fields $[v_1,[a,v_2]]$ and $[v_2,[a,v_2]]$ which is contained in $D_{n_3+1}$, \ie there exist functions $\kappa^1$ and $\kappa^2$ such that
			\begin{align*}
				\begin{aligned}
					\kappa^1[v_1,[a,v_2]]+\kappa^2[v_2,[a,v_2]]&\in D_{n_3+1}
				\end{aligned}
			\end{align*}
			with at least $\kappa^1\neq 0$ or $\kappa^2\neq 0$. At least $[v_1,[a,v_2]]\notin D_{n_3+1}$ or $[v_2,[a,v_2]]\notin D_{n_3+1}$ holds\footnote{Otherwise, we would have $\C{D_{n_3+1}}=D_{n_3}$ and thus, the condition $\C{D_{n_3+1}}\neq D_{n_3}$ of Theorem \ref{thm:1} would be violated.}. Therefore, either $[v_2,[a,v_2]]=\kappa[v_1,[a,v_2]]\Mod D_{n_3+1}$ or $[v_1,[a,v_2]]=\kappa[v_2,[a,v_2]]\Mod D_{n_3+1}$  and thus \eqref{eq:alg_1} simplifies to either
			\begin{align}\label{eq:alg_2_1}
				\begin{aligned}
					&(\alpha^1\beta_1^2+\alpha^2\beta_2^2)\left(2(\alpha^1\beta_1^1+\alpha^2\beta_2^1)+\kappa(\alpha^1\beta_1^2+\alpha^2\beta_2^2)\right)[v_1,[a,v_2]]\overset{!}{\in}D_{n_3+1}\,.
				\end{aligned}	
			\end{align}
			\vspace{-1ex}
			or
			\vspace{-1ex}
			\begin{align}\label{eq:alg_2_2}
				\begin{aligned}
					&(\alpha^1\beta_1^2+\alpha^2\beta_2^2)\left(2\kappa(\alpha^1\beta_1^1+\alpha^2\beta_2^1)+\alpha^1\beta_1^2+\alpha^2\beta_2^2\right)[v_2,[a,v_2]]\overset{!}{\in}D_{n_3+1}\,.
				\end{aligned}	
			\end{align}
			In both cases, there exist at most two independent non-trivial solutions, \ie each of the factors, which depend on $\alpha^1$ and $\alpha^2$ linearly, can vanish. In \eqref{eq:alg_2_1}, those solutions are $\alpha^1=\lambda\beta_2^2$, $\alpha^2=-\lambda\beta_1^2$ and $\alpha^1=\lambda(2\beta_2^1+\kappa\beta_2^2)$, $\alpha^2=-\lambda(2\beta_1^1+\kappa\beta_1^2)$, both with arbitrary $\lambda\neq 0$. In \eqref{eq:alg_2_2}, those solutions are $\alpha^1=\lambda\beta_2^2$, $\alpha^2=-\lambda\beta_1^2$ and $\alpha^1=\lambda(2\kappa\beta_2^1+\beta_2^2)$, $\alpha^2=-\lambda(2\kappa\beta_1^1+\beta_1^2)$, again both with arbitrary $\lambda\neq 0$. Therefore, in any case, with one of the solutions, namely $\alpha^1=\lambda\beta_2^2$ and $\alpha^2=-\lambda\beta_1^2$ with arbitrary $\lambda\neq 0$, we recover the direction of $b_1=\partial_{x_{3,1}^{n_3}}$. The second candidate $b_p$, which we obtain from the second solution (\ie $\alpha^1=\lambda(2\beta_2^1+\kappa\beta_2^2)$ and $\alpha^2=-\lambda(2\beta_1^1+\kappa\beta_1^2)$ for \eqref{eq:alg_2_1}, or  $\alpha^1=\lambda(2\kappa\beta_2^1+\beta_2^2)$ and $\alpha^2=-\lambda(2\kappa\beta_1^1+\beta_1^2)$ for \eqref{eq:alg_2_2}), may or may not be collinear with $b_1=\partial_{x_{3,1}^{n_3}}$ and the conditions of Theorem \ref{thm:1} may or may not be met with this second candidate for $b_p$.\\\\
			(If we are in Case 1, \ie if $[v_1,[a,v_2]]$ and $[v_2,[a,v_2]]$ are not collinear$\mod D_{n_3+1}$, and thus if \eqref{eq:alg_0} certainly yields only one candidate for $b_p$, can be deduced from the dimension of the distribution $H=D_{n_3+1}+[D_{n_3},D_{n_3+1}]$. To be precise, if $\Dim{H}=\Dim{D_{n_3+1}}+2$, $[v_1,[a,v_2]]$ and $[v_2,[a,v_2]]$ are not collinear$\Mod D_{n_3+1}$\footnote{We have $[v_1,[a,v_2]]=[b_1^c,b_2^c]\Mod D_{n_3+1}$ and $[v_2,[a,v_2]]=[b_2^c+gb_1^c,[a,b_2^c+gb_1^c]]$. Recall that we have $D_{n_3+1}=\Span{\partial_{x_3},b_1^c,b_2^c,[a,b_2^c+gb_1^c]}$ and $H=\Span{\partial_{x_3},b_1^c,b_2^c,[a,b_2^c+gb_1^c],[b_2^c+gb_1^c,[a,b_2^c+gb_1^c]],[b_1^c,b_2^c]}$, see \eqref{eq:H}. Thus, we actually have $H=D_{n_3+1}+\Span{[v_1,[a,v_2]],[v_2,[a,v_2]]}$. For $\Dim{H}=\Dim{D_{n_3+1}}+2$, there neither $[v_1,[a,v_2]]$ nor $[v_2,[a,v_2]]$ can already be contained in $D_{n_3+1}$, nor they can be collinear$\Mod D_{n_3+1}$.}. The Case 2 occurs if $\Dim{H}=\Dim{D_{n_3+1}}+1$. If $\Dim{D_{n_3+1}}=n-1$, $n$ being the total number of states of the system under consideration, we of course always have $\Dim{H}=\Dim{D_{n_3+1}}+1$ and in turn, we always have two candidates for $b_p$.)\\\\
			In conclusion, in any case the necessary condition \eqref{eq:alg_0} has at most two independent non-trivial solutions and thus yields at most two non-collinear candidates for the vector field $b_p$. Therefore, solving \eqref{eq:alg} to obtain candidates for $b_p$, needed for applying Theorem \ref{thm:1}, we obtain at most two non-collinear candidates, and if the system is indeed static feedback equivalent to \eqref{eq:triangularFormDifferentLengths}, then at least for one of these candidates, the conditions of Theorem \ref{thm:1} are met.
\end{document}